\makeatletter
\RequirePackage{keyval}
\define@key{Gm}{twosideshift}{%
	\ifdim#1=0pt%
	\else
	\PackageWarning{geometry}{Option `twosideshift` is gone since version 5}%
	\fi
}
\makeatother
\providecommand*{\StandardLayout}{}

\documentclass[envcountsame]{my-mmnp}

\usepackage[tbtags]{amsmath}
\usepackage{amssymb}
\usepackage{graphicx}
\usepackage[latin1]{inputenc}

\usepackage{subfig}

\idline{Vol. 12, No. 3}{51}
\doi{10.1051/mmnp/201712305}

\numberwithin{equation}{section}

\begin{document}


\title{Optimal Spraying in Biological Control of Pests}

\author{C. J. Silva \inst{1} 
\sep D. F. M. Torres \inst{1}\thanks{\email{delfim@ua.pt}} 
\sep E. Venturino \inst{2}}


\vspace{0.5cm}

\institute{\inst{1} Center for Research and Development in Mathematics and Applications (CIDMA),\\
Department of Mathematics, University of Aveiro, 3810--193 Aveiro, Portugal.\\
\inst{2} Dipartimento di Matematica ``Giuseppe Peano'',
Universit\`{a} di Torino, Italy.}


\abstract{We use optimal control theory with the purpose of finding 
the best spraying policy with the aim of at least to minimize
and possibly to eradicate the number of parasites, i.e.,
the prey for the spiders living in an agroecosystems.
Two different optimal control problems are posed and solved,
and their implications discussed.}


	



\titlerunning{Optimal spraying in biological control of pests}

\maketitle


\section{Introduction}

Fairly recent researches have highlighted the role
that generalist predators, i.e., predators that feed 
on several species, have on pests in agriculture,
so that their presence in agroecosystems should be fostered \cite{RL}.
For instance, it is well known that toward the end of the 
nineteenth century, the European vineyards \emph{Vitis vinifera}
suffered from the accidental introduction
of the plant louse \emph{phylloxera} \cite{Pinney:HW}.
Later vineyard agroecosystems in every continent became affected by
this spread of the pest. This fact caused relevant 
economic damages worldwide \cite{Winemaking}.
The crisis was overcome by the use of native American grapevines,
grafting the European vines onto Eastern American roots, such as
\emph{Vitis labrusca}, which are resistant to \emph{phylloxera} \cite{Ox:Comp:Wine}.
The role of spiders to combat these and other parasites affecting vines,
such as black measles, little-leafs, nematodes and red ticks, and
infesting fruit orchards, e.g., grapevine beetles, grape-berry moths, 
climbing cutworms, black rot and mildew, has been recently elucidated 
in the literature \cite{CD98,IBBB,NM,W}.
The findings of these researches indicate that
vineyards contain a very diversified and abundant spider community,
which depends heavily on the availability of the grass lying in between the vine rows
and possibly nearby bushes or small woods,
and as a whole helps in keeping down the parasites population \cite{RB}.
However, spiders have different hunting strategies
and locations, for which different insects are most commonly hunted by the various
spider species \cite{MC}.
Spiders can be distinguished in two very broad sets, the
web-builder species, which are in general stantial, except for a single ``flight''
due to the action of the wind when they are young.
Then they build the web and wait for the prey. The second set is composed by
wanderer species, like wolf spiders, that crawl on the ground in search for food.
For instance, the nocturnal wandering spiders \emph{Clubiona brevipes}, 
\emph{C. corticalis} and \emph{C. leucaspis} (Clubionidae)
hunt non-flying Aphids and larvae of Lepidoptera, while the diurnal
wandering species \emph{Ballus depressus} (Salticidae) are effective against
Cicadellidae; adults and larvae of Hymenoptera and Lepidoptera are prey 
of the ambush species \emph{Philodromus aureolus}
(Philodromidae) and \emph{Diaea dorsata} (Thomisidae) \cite{IBB06,Ma}.
The ultimate result is that we can consider spiders as specialists,
i.e., feeding essentially only on one species,
thus implying that parasites can be controlled only by ensuring the persistence
of a wide variety of diverse spider species within the agroecosystem.
Mathematical models for the understanding of these complex systems have been
proposed and analysed in the past years: in \cite{VIBITB} a space-free
description of the dispersal phenomenon of web-builder spiders is presented;
wanderer spiders are instead considered in \cite{VIBCB} and its extensions
\cite{CIVa,CIVb}.

Another very common strategy to fight parasites is represented
by the use of pesticides, although it is now widely
recognized that they have high economical costs and negative impact
in terms of both environment damages
as well as human health. Besides, pests acquire resistance to the chemicals,
which may lead to unwanted large parasites outbreaks \cite{RD-B}.

The action of parasites has been also considered within the modeling framework
of previous researches, but always considering very simple control strategies,
as the focus was mainly on extracting the relationships between the various populations
living in the environment, while the human role within the ecosystem model was considered
marginal. In fact, either an instantaneous effect of the spray, \cite{VIBITB},
or an exponentially decay of the poison, \cite{CIVb,VIBCB}, were considered, but in both
cases the control was taken to be a constant function.

Optimal control theory has a long history of being applied to problems in biomedicine
(see, e.g., \cite{Eisen1979,Ledzewicz2002,Lenhart2007} and references cited therein)
and recently to epidemic models for infectious diseases (see, e.g.,
\cite{Behncke2000,Gaff2009,Sofia2010,Sofia2013,PaulaSilvaTorresT2014,%
SilvaTorresNACO,SilvaTorresMBS2013,Silva:Torres:DCDS-A:2015} for population compartmental models, \cite{HattafHIV2012,HattafHIVdelay2012} for HIV cellular-level models, 
and \cite{IbrahimIJDC2016} for numerical methods applied to optimal control problems). 
However, to our knowledge, little attention has been given to agroecosystem models.

In this paper, we apply optimal control theory to the activities related to an agroecosystem,
and study specifically the possible spraying policies of the fruit orchards.
Two optimal control problems are identified and solved in this context,
providing optimal strategies for the realistic implementation of pest eradication:
the complete elimination of the pests, the minimum time to achieve the previous goal
and the minimization of the parasites in a given time frame.

The text is organized as follows. In Section~\ref{sec:model}, we recall the model
and extend it as a general control system. In Section~\ref{sec:2}, two optimal
control problems are proposed and solved: subsection~\ref{sub:sec:ocp1}
addresses the pest eradication problem, subsection~\ref{sec:ocp2} contains 
the time optimal control problem. We end with Section~\ref{sec:disc} of discussion
and Section~\ref{sec:conc:fw} of conclusions and future work.


\section{Mathematical model with control function}
\label{sec:model}

We consider the mathematical model from \cite{Ezio:spiders:JNAIAM:2008},
which considers interactions between three populations: insects living
in open fields and woods at time $t$ ($f(t))$; parasites living
in the vineyards ($v(t)$); and wanderer spiders living in the
whole environment ($s(t)$). Both insects $f(t)$ and parasites $v(t)$
can constitute prey for the spiders $s(t)$. The time unit is the day and
populations are counted in numbers.

Let $W$ and $V$ denote the carrying capacity for insects in the woods and vineyards,
respectively. Given the situation of highly exploited lands, it is assumed
that the vineyard carrying capacity much exceeds the one of nearby
green areas, so that $V \gg W$.

We assume that insect and parasite populations ($f$ and $v$) are the only
food source for spiders $w$, so that when the former are absent
the latter will die exponentially fast. Predation
is accounted for by mass action law terms with suitable signs.
Captured prey contribute to spider reproduction, however not the whole prey
is consumed but only a fraction of it; this is expressed 
via the efficiency constant $0 < k < 1$.

Insects and parasites natural birth rates are denoted by $r$ and $e$, respectively.
The spiders' natural mortality is denoted by the parameter $a$. The hunting rates
of spiders on vineyard and woods are denoted by $b$ and $c$, respectively. The parameter
$h$ represents the intensity of spraying, i.e., the amount of poison released, $K$ models
its unwanted killing effect on spiders, $q$ represents the fraction of insecticide that
lands on target, i.e., within the vineyards, and $1-q$ the corresponding fraction that due
to the action of the wind, water and possibly other causes, is dispersed in the woods.
All the parameters are nonnegative constants unless otherwise specified.

We rewrite the most general form of the model in \cite{Ezio:spiders:JNAIAM:2008}
including human spraying of the fields, by replacing 
its particular controls by a general function.
Specifically, we add to the dynamic system a control variable, $u$, which is
a function of time and not a very specific control as assumed in the
previous literature \cite{Ezio:spiders:JNAIAM:2008,CIVb,VIBCB,VIBITB}.
It represents the human intervention \emph{via} spraying with insecticides, as follows:
\begin{equation}
\label{model:controls}
\begin{cases}
\dot{f}(t) = r f(t)\left(1 - \frac{f(t)}{W}\right) - c s(t) f(t) - h(1-q) u(t),\\[0.2 cm]
\dot{s}(t) = s(t)\left(-a + k b v(t) + k c f(t)\right) - h K q u(t),\\[0.2 cm]
\dot{v}(t) = e v(t)\left(1-\frac{v(t)}{V}\right)- b s(t)v(t) - h q u(t).
\end{cases}
\end{equation}


\section{Problem formulation}
\label{sec:2}

In this section, we formulate two optimal control problems where the control system
is given by \eqref{model:controls}. These problems are subject to initial conditions
\begin{equation}
\label{init:cond}
f(0) = f_0   \, , \quad s(0) = s_0  \, ,
\quad   v(0) =  v_0  \, ,
\quad \text{with} \, \, f_0, s_0, v_0 \geq 0 \, .
\end{equation}
Note that for $u \equiv 0$ the control system \eqref{model:controls} reduces to
the model for the ecosystem without human intervention.
The set $\Omega$ of admissible control functions is defined as
\begin{equation*}
\Omega := \left\{ u(\cdot) \in L^{\infty}(0, T) \,
| \,  0 \leq u(t) \leq 1 ,  \, \forall \, t \in [0, T] \, \right\} .
\end{equation*}
The optimal control problems are the following:
\begin{enumerate}
\item minimize the number of parasites ($v(t)$) with human intervention (use of insecticides),
taking or not into account the cost of insecticides (see subsection~\ref{sub:sec:ocp1});
	
\item minimize the time $T$ for which the parasites of
the vineyard will be eradicated, that is, $v(T) = 0$ with $f(T)$ and $s(T)$ free and $s(T) > 0$
(see subsection~\ref{sec:ocp2}).
	
\end{enumerate}
For the two problems, we compare the results with human intervention (with control, that is, $u \ne 0$)
to the case where there is no human intervention (no control/no insecticides, that is, $u \equiv 0$).


\subsection{Minimizing the number of parasites}
\label{sub:sec:ocp1}

Consider the following objective functional defined by a sum of two terms:
\begin{equation}
\label{costfunction}
\mathcal{J}(u(\cdot)) = \int_0^{T} \left[ v(t) + \frac{\xi}{2}u^2(t) \right] dt.
\end{equation}
The first term in our cost functional tells us that we want to minimize the number of pests.
The constant $\xi$ in the second term is a measure of the relative cost of the intervention
associated to the control $u$. One applies control measures that are associated with some implementation
costs that we also intend to minimize. By considering the cost with the control in a quadratic form,
we are being consistent with previous works in the literature
(see, e.g., \cite{PaulaSilvaTorresT2014,SilvaTorresMBS2013}).
Moreover, a quadratic structure in the control has mathematical advantages.
Roughly speaking, it implies that the Hamiltonian attains its minimum over
the control set at a unique point (given by \eqref{optcontrols}).
The optimal control problem consists of determining
the triple $\left(f^*(\cdot), s^*(\cdot), v^*(\cdot)\right)$,
associated to an admissible control
$u^*(\cdot) \in \Omega$ on the time interval $[0, T]$,
satisfying \eqref{model:controls},
the initial conditions \eqref{init:cond}, and
minimizing the cost functional \eqref{costfunction}, that is,
\begin{equation}
\label{mincostfunct}
\mathcal{J}(u^*(\cdot))
= \min_{\Omega} \mathcal{J}(u(\cdot)) \, .
\end{equation}

The existence of an optimal control $u^*$ and associated $(f^*, s^*, v^*)$
comes from the convexity of the integrand of the cost functional \eqref{costfunction} with
respect to the control $u$ and the Lipschitz property of the state system with respect
to state variables $(f, s, v)$ (see, e.g.,  \cite{Cesari_1983,Fleming_Rishel_1975}
for existence results of optimal solutions).

According to the Pontryagin Maximum Principle \cite{Pontryagin_et_all_1962},
if $u^*(\cdot) \in \Omega$ is optimal for the problem \eqref{model:controls},
\eqref{mincostfunct} with the initial conditions given by \eqref{init:cond}
and fixed final time $T$, then there exists a nontrivial absolutely
continuous mapping $\lambda : [0, T] \to \mathbb{R}^3$,
$\lambda(t) = \left(\lambda_1(t), \lambda_2(t), \lambda_3(t)\right)$,
called \emph{adjoint vector}, such that
\begin{equation*}
\dot{f} = \frac{\partial H}{\partial \lambda_1} \, , \quad
\dot{s}= \frac{\partial H}{\partial \lambda_2} \, , \quad
\dot{v}= \frac{\partial H}{\partial \lambda_3}
\end{equation*}
and
\begin{equation}
\label{adjsystemPMP}
\dot{\lambda}_1 = -\frac{\partial H}{\partial f} \, , \quad
\dot{\lambda}_2 = -\frac{\partial H}{\partial s} \, , \quad
\dot{\lambda}_3 = -\frac{\partial H}{\partial v} \, ,
\end{equation}
where the function $H$ defined by
\begin{multline*}
H = H(f, s, v, u, \lambda_1, \lambda_2, \lambda_3)
= v + \frac{\xi}{2}u^2 +\lambda_1 \left( r f\left(1
- \frac{f}{W}\right) - c s f - h(1-q) u \right)\\
+\lambda_2 \left( s\left(-a + k b v + k c f\right) - h K q u \right)
+\lambda_3 \left( e v\left(1-\frac{v}{V}\right)- b s v - h q u \right)
\end{multline*}
is called the \emph{Hamiltonian}, and the minimality condition
\begin{equation}
\label{maxcondPMP}
H(f^*(t), s^*(t), v^*(t), u^*(t), \lambda(t))
= \min_{0 \leq u \leq 1} H(f^*(t), s^*(t), v^*(t), u, \lambda(t))
\end{equation}
holds almost everywhere on $[0, T]$. Moreover, one has
the transversality conditions $\lambda_i(T) = 0$, $i =1,\ldots, 3$.

\begin{thrm}
\label{the:thm}
The problem \eqref{model:controls}, \eqref{mincostfunct} with fixed initial conditions
\eqref{init:cond} and fixed small final time $T$, admits a unique optimal solution
$(f^*(\cdot), s^*(\cdot), v^*(\cdot))$ associated with an optimal control
$u^*(\cdot)$ on $[0, T]$. Moreover, there exist adjoint functions,
$\lambda_1(\cdot)$, $\lambda_2(\cdot)$ and $\lambda_3(\cdot)$, such that
\begin{equation}
\label{adjoint_function}
\begin{cases}
\dot{\lambda}_1(t) = -\lambda_1(t) \left( r \left( 1-{\frac {f^*(t)}{W}} \right)
-{\frac {rf^*(t)}{W}}-c s^*(t) \right) - \lambda_2(t) \,s^*(t) k c,\\
\dot{\lambda}_2(t) = \lambda_1(t)\,cf^*(t) - \lambda_2(t) \,
\left( -a + k b v^*(t)+k c f^*(t) \right) + \lambda_3(t) \,b v^*(t),\\
\dot{\lambda}_3(t) = -1- \lambda_2(t) \,s^*(t) k b- \lambda_3(t)
\, \left( e \left( 1-{\frac {v^*(t)}{V}} \right) -{\frac {e v^*(t)}{V}}-b s^*(t) \right),
\end{cases}
\end{equation}
with transversality conditions $\lambda_i(T)=0$, $i=1, \ldots, 3$. Furthermore,
\begin{equation}
\label{optcontrols}
u^*(t) = \min \left\{ \max \left\{0, {\frac {h \left( \lambda_1(t) - \lambda_1(t) \,q
+ \lambda_2(t) \,K q+ \lambda_3(t) \,q \right) }{\xi}}  \right\}, 1 \right\} \, .
\end{equation}
\end{thrm}

\begin{proof}
Existence of an optimal solution $(f^*, s^*, v^*)$
associated to an optimal control $u^*$ comes from
the convexity of the integrand of the cost functional $J$ with respect
to the control $u$ and the Lipschitz property of the state system
with respect to the state variables $\left(f, s, v \right)$
(see, \textrm{e.g.}, \cite{Cesari_1983,Fleming_Rishel_1975}).
System \eqref{adjoint_function} is derived from the Pontryagin maximum principle
(see \eqref{adjsystemPMP}, \cite{Pontryagin_et_all_1962})
and the optimal controls \eqref{optcontrols} come from the minimization condition \eqref{maxcondPMP}.
For small final time $T$, the optimal control pair given by \eqref{optcontrols}
is unique due to the boundedness of the state and adjoint functions and the Lipschitz property
of systems \eqref{model:controls} and \eqref{adjoint_function}
(see \cite{SLenhart_2002} and references cited therein).
\end{proof}

In this article the numerical results were obtained using
the PROPT Matlab Optimal Control Software \cite{PROPT}. 
We considered the initial conditions
\begin{equation}
\label{eq:init:cond:p1}
f(0) = 3.1   \, , \quad s(0) = 3.7  \, , \quad   v(0) =  2.2   \, ,
\end{equation}
and the parameter values from Table~\ref{parameters}.
\begin{table}[!htb]
\centering
\begin{tabular}{|l | l | l |}
\hline
{\small{Symbol}} & {\small{Description}}  & {\small{Value}} \\
\hline
{\small{$r$}} & {\small{insects's natural birth rate}} & {\small{1}} \\
{\small{$e$}} & {\small{parasites's natural birth rate}} & {\small{2.5}} \\
{\small{$a$}} & {\small{spider's natural mortality}} & {\small{3.1}} \\
{\small{$b$}} & {\small{hunting rate of spiders on vineyard insects}} & {\small{1.2}} \\
{\small{$c$}} & {\small{hunting rate of spiders on wood insects}} & {\small{0.2}} \\
{\small{$h$}} & {\small{intensity of spraying}} & {\small{0.7}} \\
{\small{$q$}} & {\small{fraction of insecticide that lands on vineyards}} & {\small{0.9}}\\
{\small{$k$}} & {\small{``efficiency'' constant at which prey biomass is turned into new spiders}} & {\small{1}}\\
{\small{$W$}} & {\small{carrying capacity for insects in the woods}} & {\small{5}} \\
{\small{$V$}} & {\small{carrying capacity for insects in the vineyards}} & {\small{1000}} \\
{\small{$K$}} & {\small{unwanted killing effect of spraying on spiders}} & {\small{0.01}} \\
\hline
\end{tabular}
\caption{Model parameters: the reference values are taken from \cite{Ezio:spiders:JNAIAM:2008}.}
\label{parameters}
\end{table}

\begin{figure}[!htb]
\centering
\subfloat[$f(t)$]{\label{fig:preyswood1}
\includegraphics[scale=0.33]{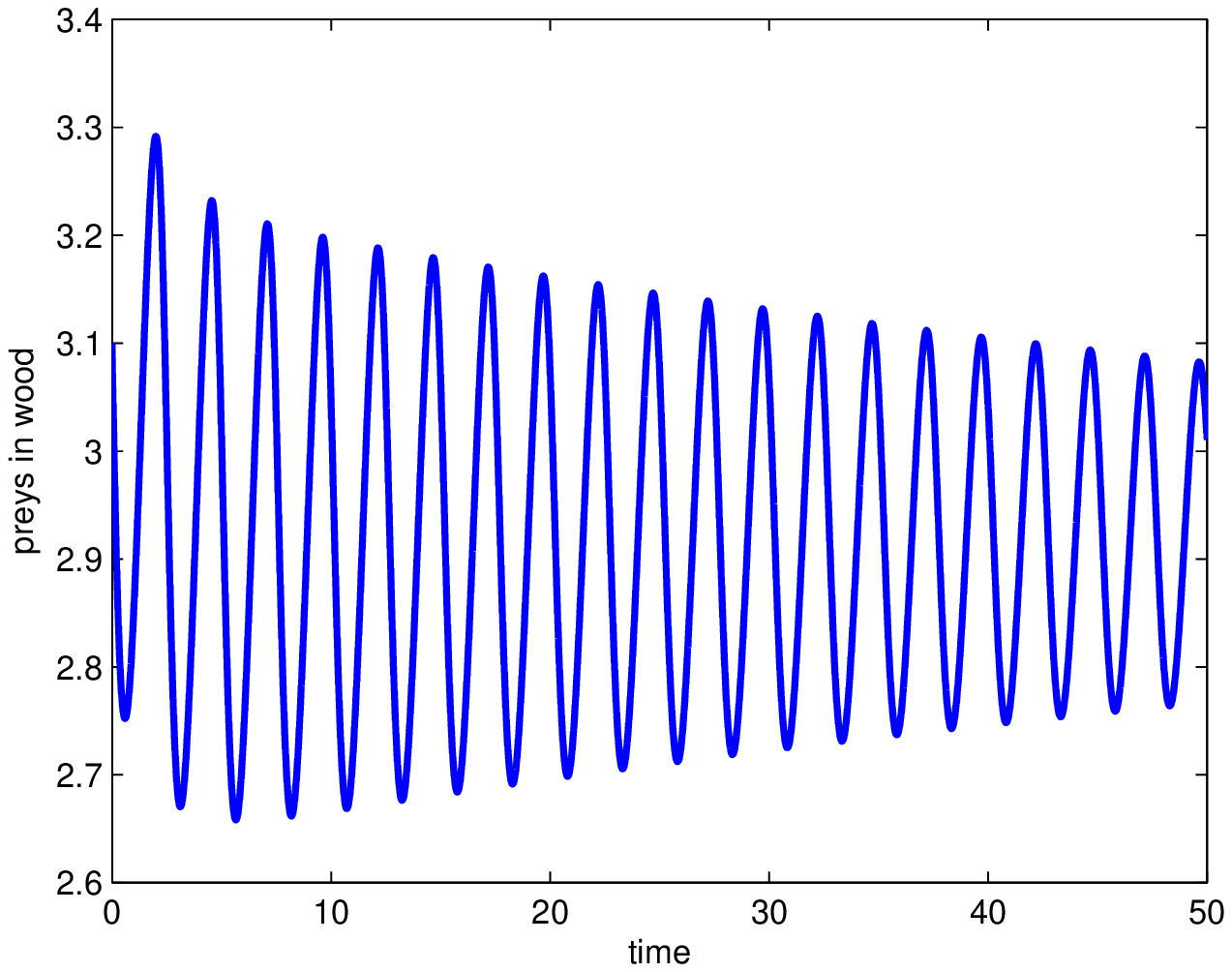}}
\subfloat[$s(t)$]{\label{fig:spiders1}
\includegraphics[scale=0.33]{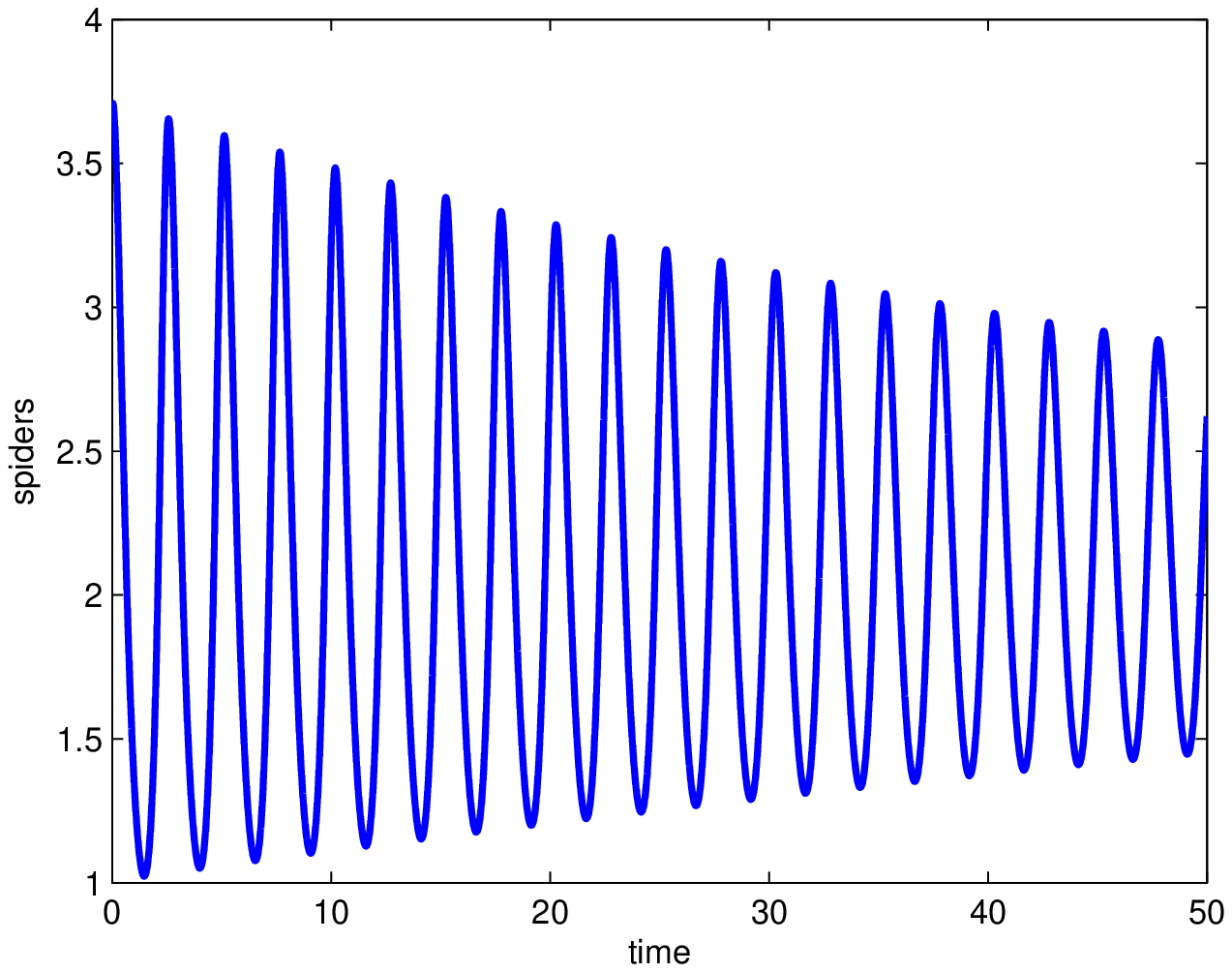}}
\subfloat[$v(t)$]{\label{fig:preyswine1}
\includegraphics[scale=0.33]{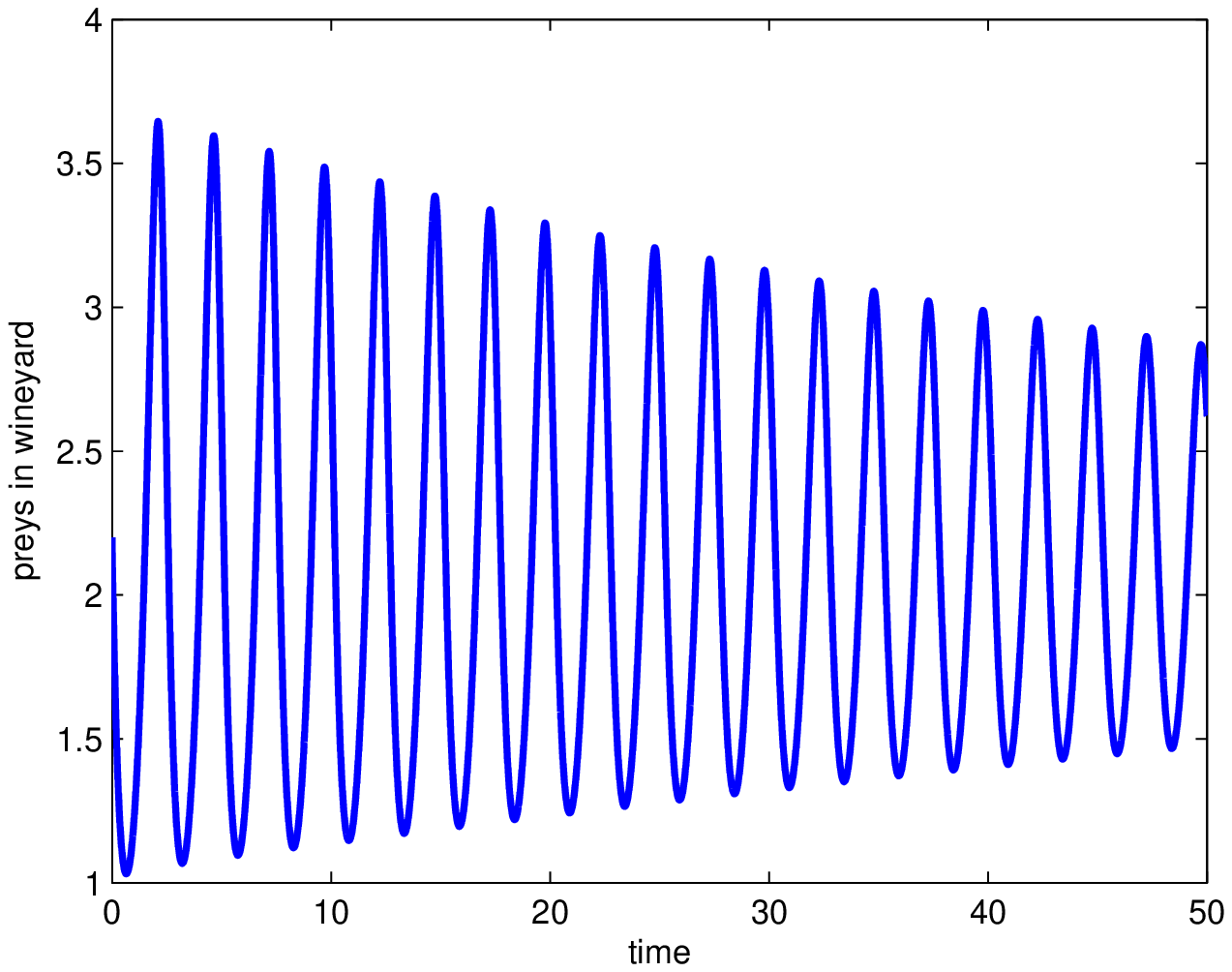}}
\caption{Behavior of \eqref{model:controls}--\eqref{init:cond}
with data given in \eqref{eq:init:cond:p1} and Table~\ref{parameters}
without spraying ($u \equiv 0$). Here and in all subsequent figures, we plot left
to right the three populations: (a) of insects in the woods
$f$, (b) spiders $s$ and (c) pests in the vineyards $v$.}
\label{fig:no:spray}
\end{figure}

If we consider the model without controls, i.e.,
\eqref{model:controls} with $u \equiv 0$,
then the behavior for $t \in [0, 50]$
is shown in Figure~\ref{fig:no:spray}.

Next, we consider human intervention, that is, the situation when
spraying with insecticides is considered ($u \ne 0$).


\subsubsection{Insecticides cost not accounted for}
\label{sec:xi0}

\begin{figure}[!htb]
\centering
\subfloat[$f^*(t)$]{\label{fig:preyswood:xi0:tf50}
\includegraphics[scale=0.33]{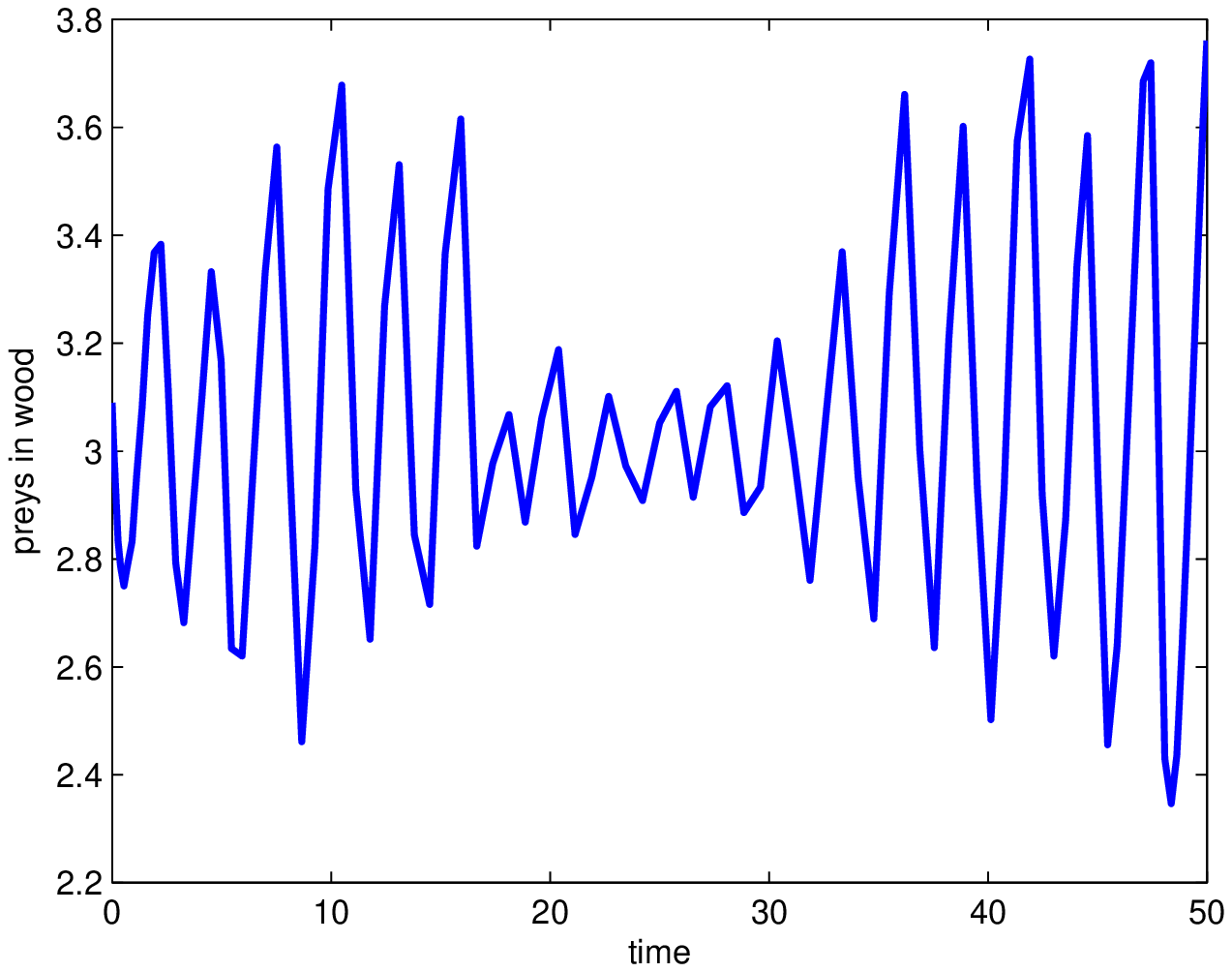}}
\subfloat[$s^*(t)$]{\label{fig:spiders:xi0:tf50}
\includegraphics[scale=0.33]{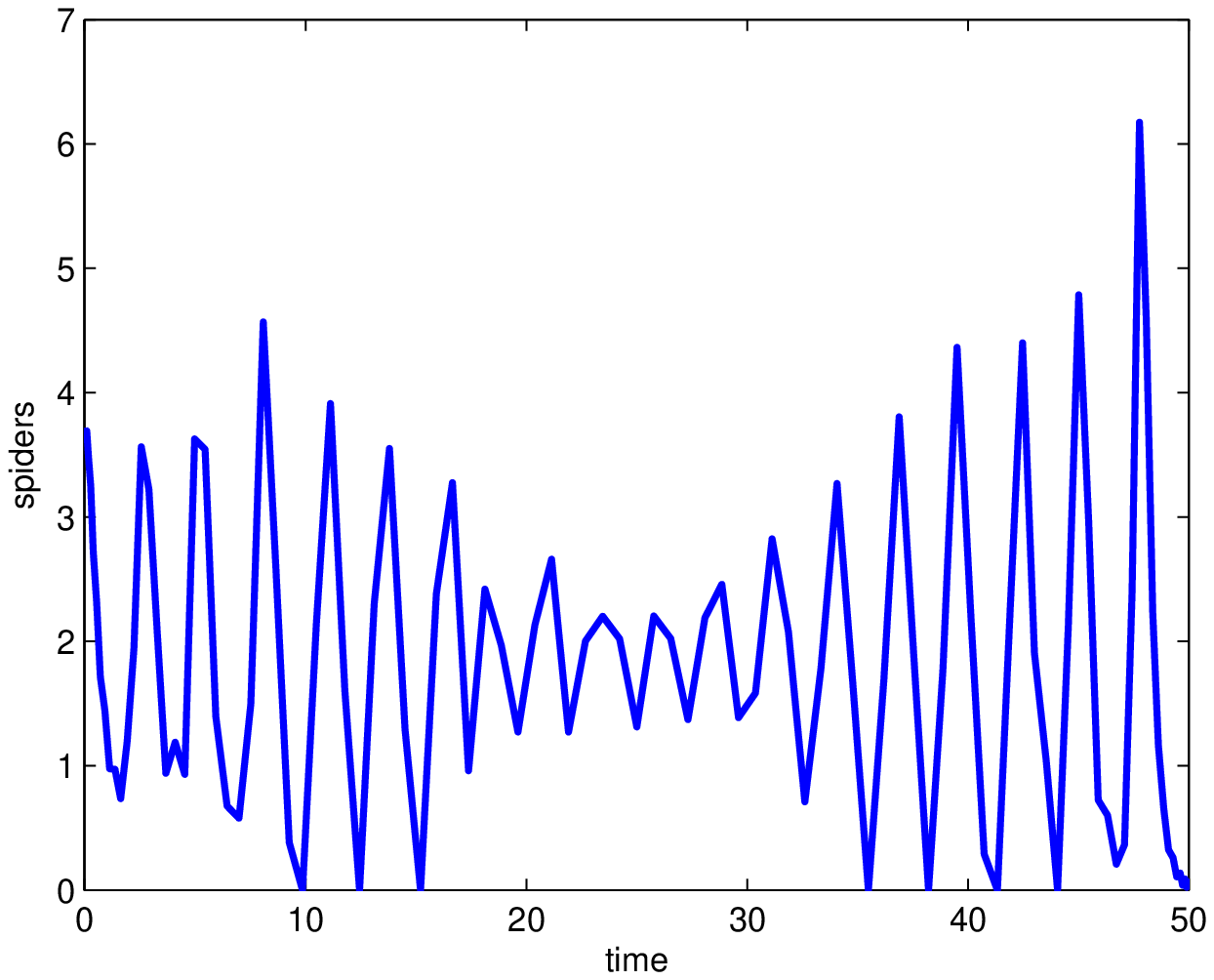}}
\subfloat[$v^*(t)$]{\label{fig:preyswine:xi:tf50}
\includegraphics[scale=0.33]{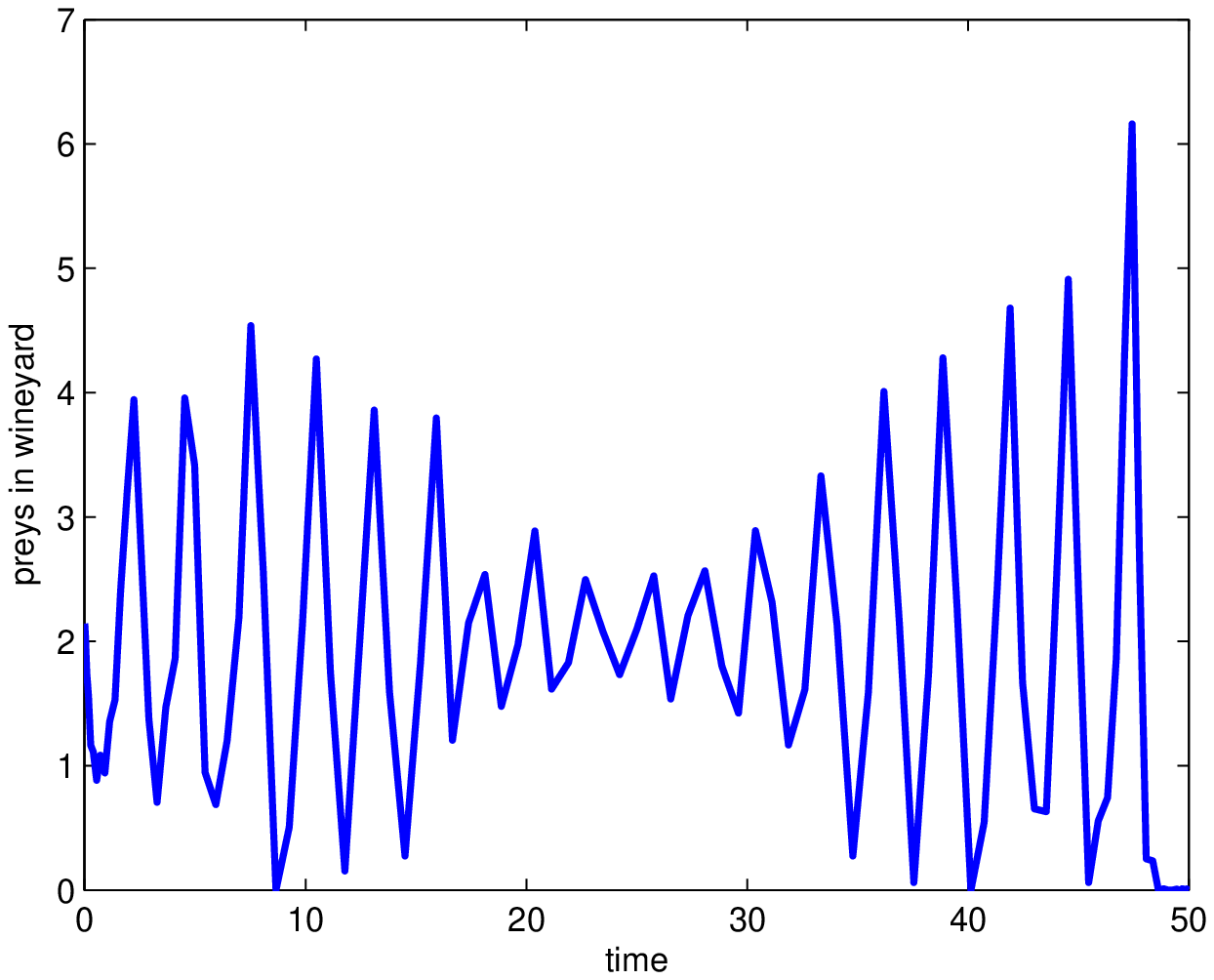}}
\caption{Solution of the optimal control problem
\eqref{model:controls}--\eqref{eq:init:cond:p1}
and parameter values from Table~\ref{parameters}
with $\xi = 0$ and $T = 50$ (Section~\ref{sec:xi0}).
The corresponding optimal control $u^*(t)$
is shown in Figure~\ref{fig:control:xi0}.}
\label{fig:spray:xi0}
\end{figure}

\begin{figure}[!htb]
\centering
\includegraphics[scale=0.49]{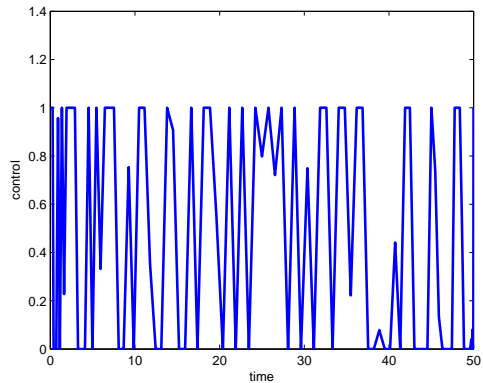}
\caption{The optimal spraying $u^*(t)$ for problem
\eqref{model:controls}--\eqref{eq:init:cond:p1}
and parameter values from Table~\ref{parameters}
with $\xi = 0$ and $T = 50$ (Section~\ref{sec:xi0}).}
\label{fig:control:xi0}
\end{figure}

\begin{figure}[!htb]
\centering
\subfloat[$f^*(t)$]{\label{fig:preyswood:xi0:tf150}
\includegraphics[scale=0.33]{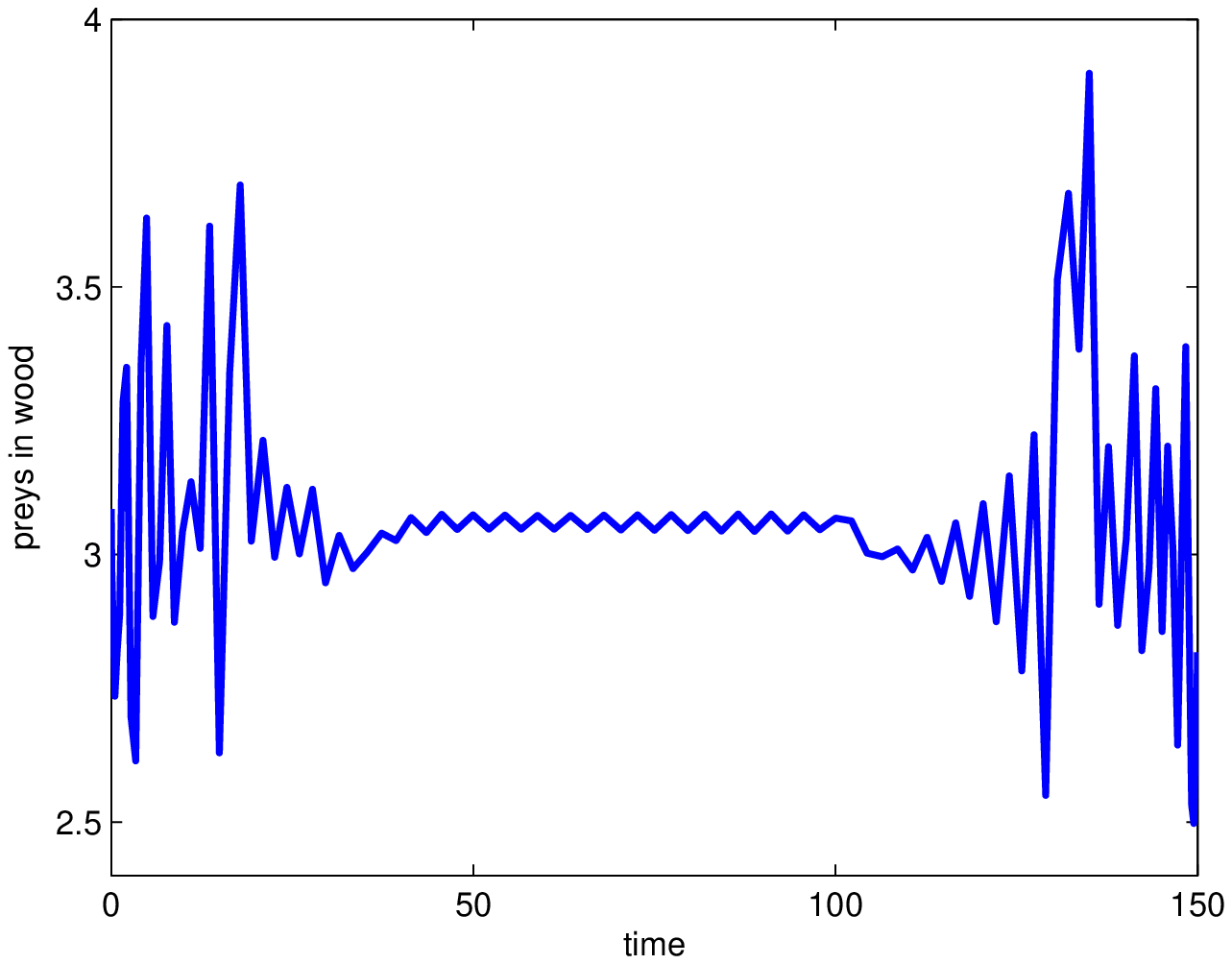}}
\subfloat[$s^*(t)$]{\label{fig:spiders:xi0:tf150}
\includegraphics[scale=0.33]{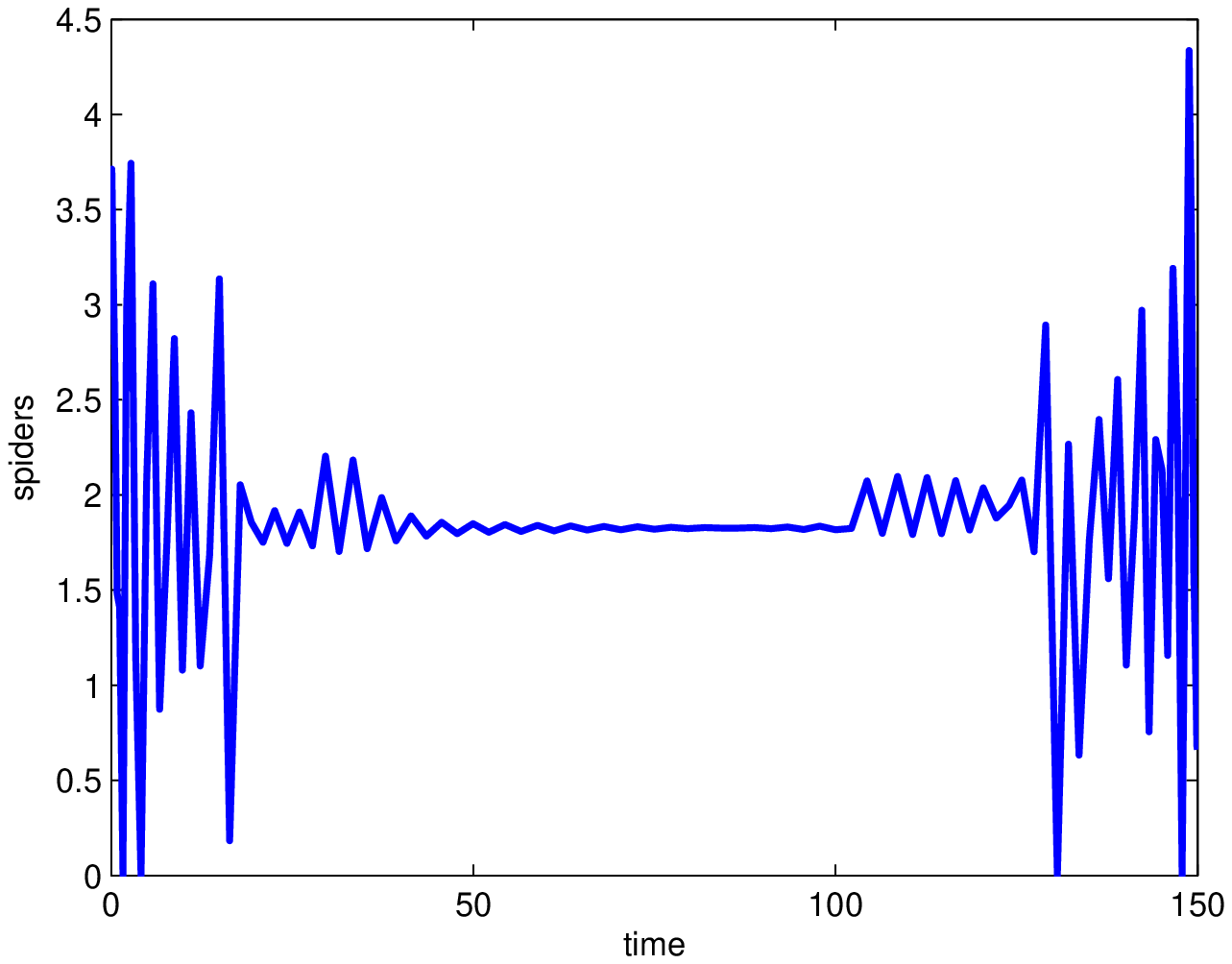}}
\subfloat[$v^*(t)$]{\label{fig:preyswine:xi:tf150}
\includegraphics[scale=0.33]{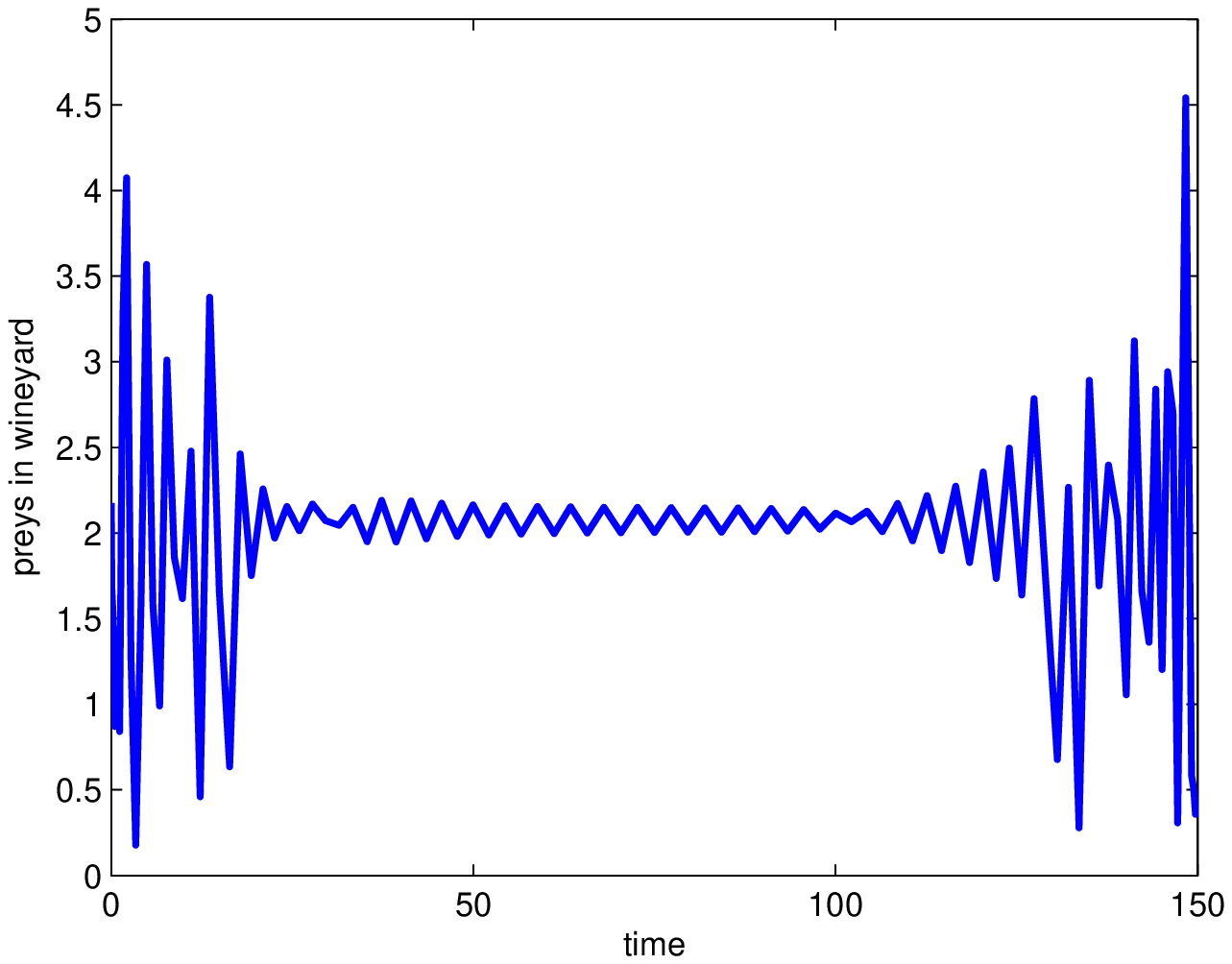}}
\caption{Solution to the optimal control problem
\eqref{model:controls}--\eqref{eq:init:cond:p1}
and parameter values from Table~\ref{parameters}
with $\xi = 0$ and $T = 150$ (Section~\ref{sec:xi0}).
The corresponding optimal control $u^*(t)$ is shown
in Figure~\ref{fig:control:xi0:tf150}.}
\label{fig:spray:xi0:tf150}
\end{figure}

\begin{figure}[!htb]
\centering
\includegraphics[scale=0.49]{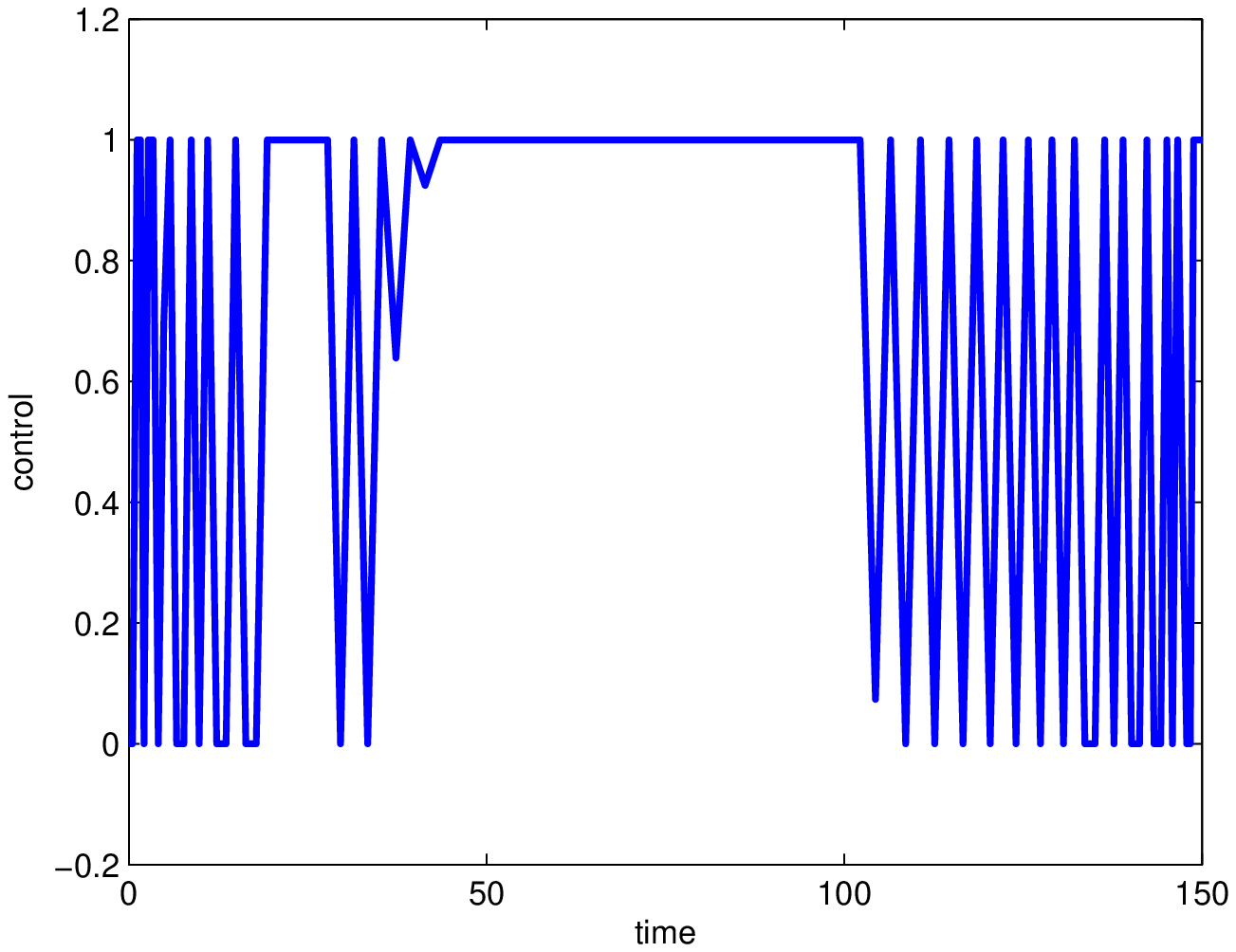}
\caption{The optimal spraying $u^*(t)$ for problem 
\eqref{model:controls}--\eqref{eq:init:cond:p1}
and parameter values from Table~\ref{parameters}
with $\xi = 0$ and $T = 150$ (Section~\ref{sec:xi0}).}
\label{fig:control:xi0:tf150}
\end{figure}

We begin by considering $\xi = 0$, that is, the goal is simply to minimize
$$
\int_0^T v(t) dt.
$$
The solution to this optimal control problem is shown in Figures~\ref{fig:spray:xi0}
and \ref{fig:control:xi0} for $T = 50$; and in Figures~\ref{fig:spray:xi0:tf150} and
\ref{fig:control:xi0:tf150} for $T = 150$.
When we do not consider the costs associated to the control (the spray),
the number of pests in the vineyard $v(t)$ can attain the value
six for $T$ near 50 days (see Figure~\ref{fig:spray:xi0}).
However, between days 20 and 30 the number $v(t)$
is lower than the one observed in Figure~\ref{fig:no:spray},
when no control is applied. Since there is no cost associated with the spray, the
optimal control attains the upper bound often (see Figure~\ref{fig:control:xi0}).
In Figures~\ref{fig:spray:xi0:tf150} and \ref{fig:control:xi0:tf150} we consider
$\xi = 0$ and $T = 150$. We observe that the number of preys in the vineyard
attains a value around 2 during almost 100 days, which is lower than the one
observed when no spray is applied (compare Figures~\ref{fig:no:spray}
and \ref{fig:spray:xi0:tf150}).
However, this is associated with a control that takes constantly
the maximum value for more than 50 days (see Figure~\ref{fig:control:xi0:tf150}).


\subsubsection{Inclusion of the cost of insecticides}
\label{sec:pb1:xineq0}

\begin{figure}[!htb]
\centering
\subfloat[$f^*(t)$]{\label{fig:preyswood:op1}
\includegraphics[scale=0.33]{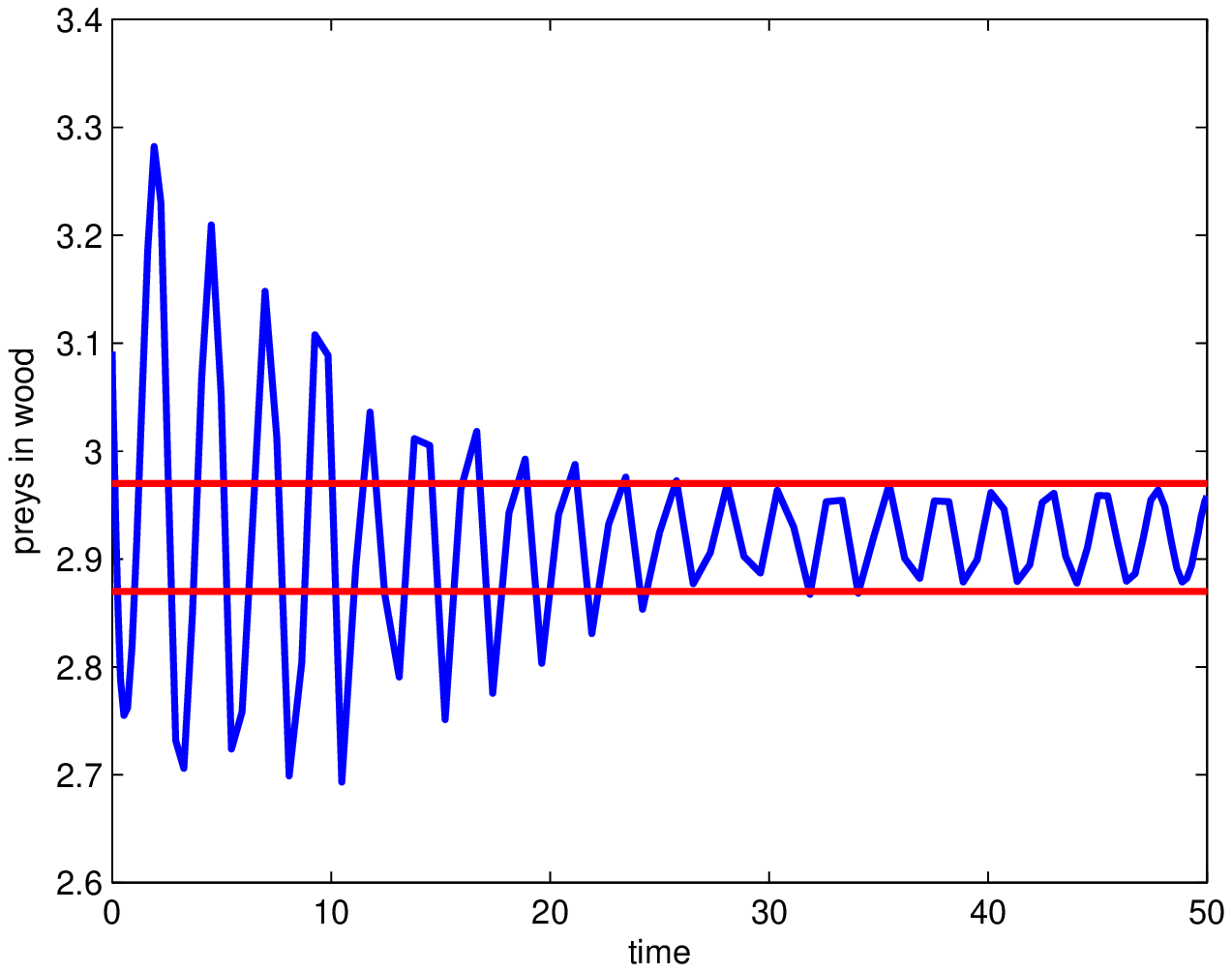}}
\subfloat[$s^*(t)$]{\label{fig:spiders:op1}
\includegraphics[scale=0.33]{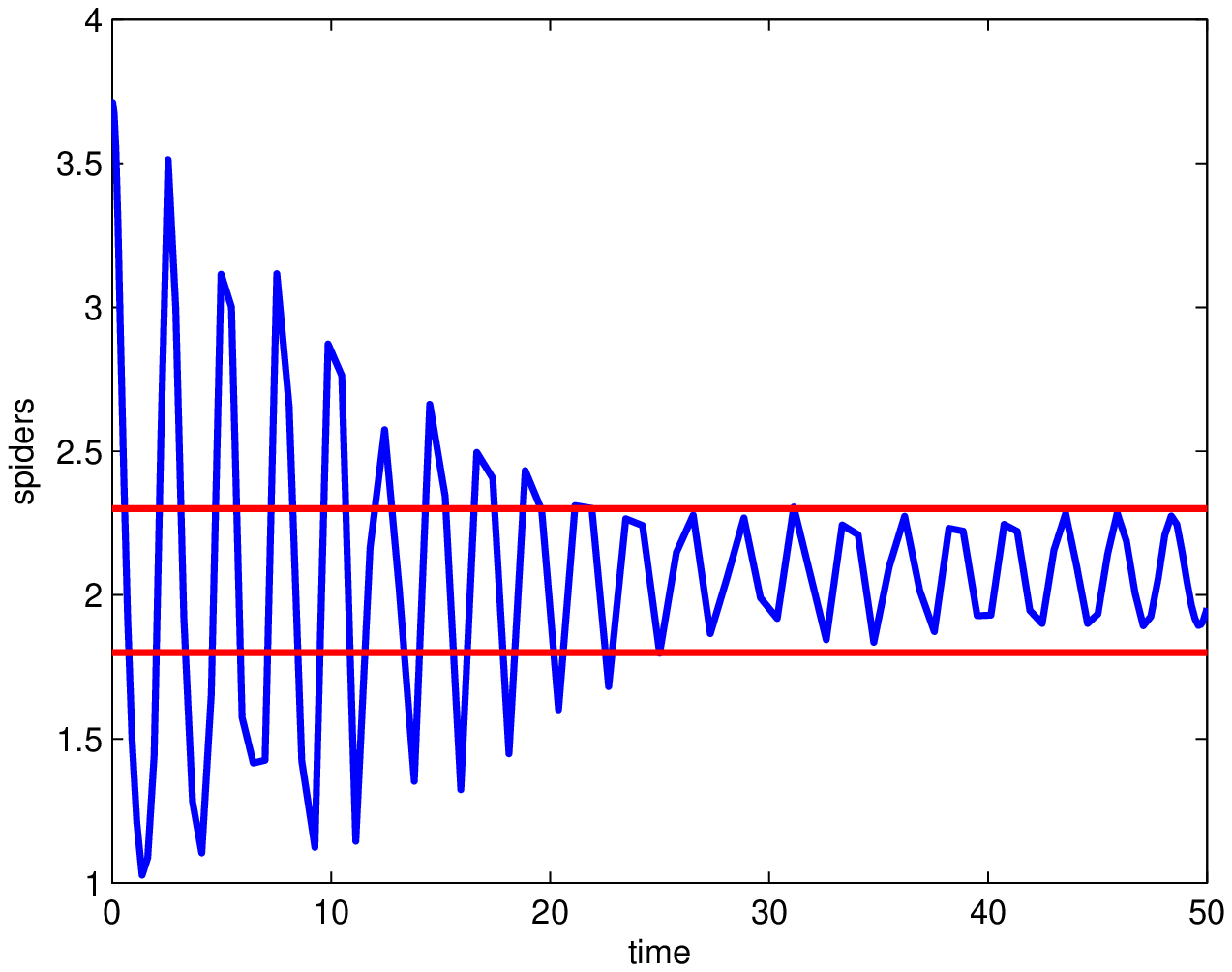}}
\subfloat[$v^*(t)$]{\label{fig:preyswine:op1}
\includegraphics[scale=0.33]{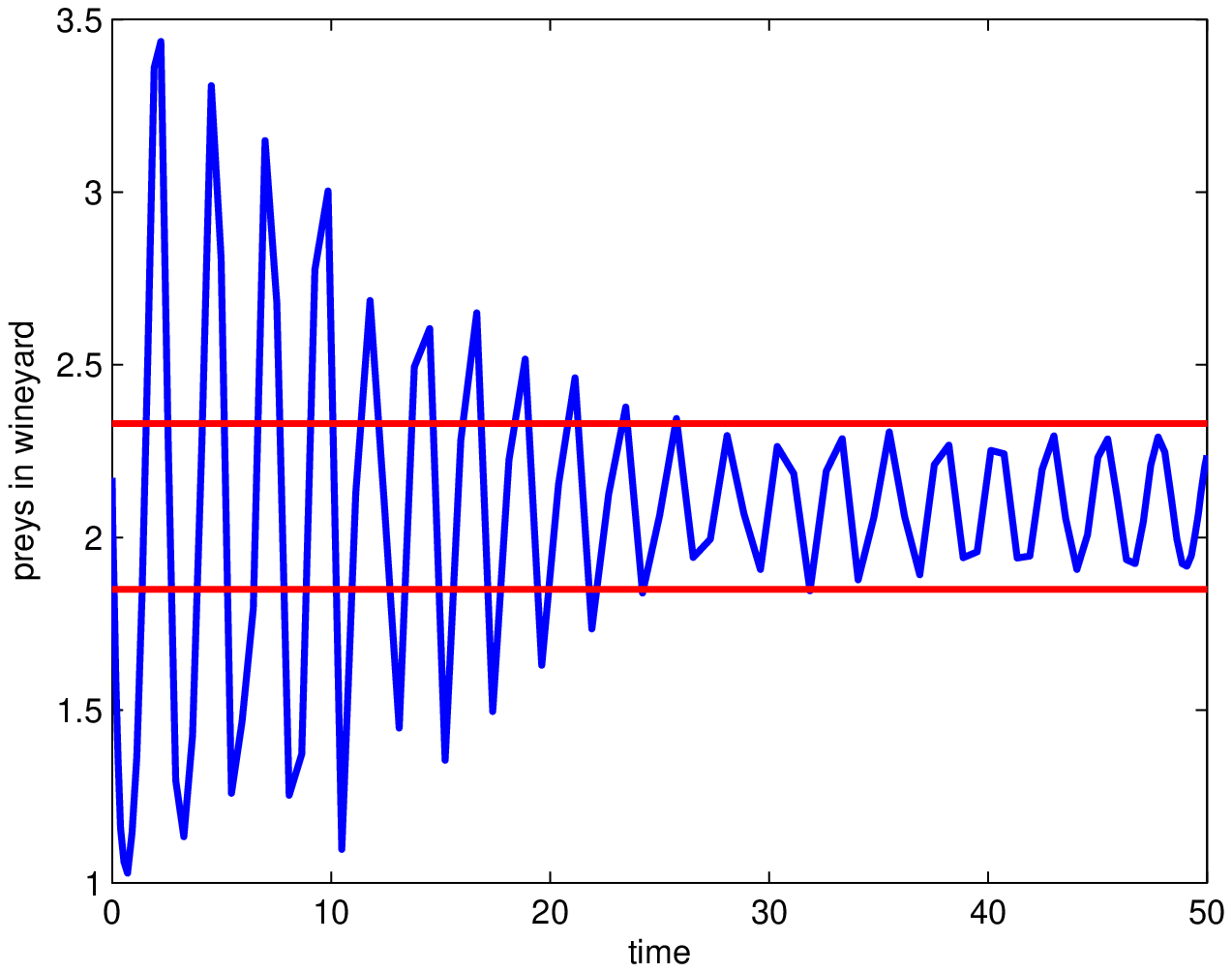}}
\caption{Solution of the optimal control problem
\eqref{model:controls}--\eqref{eq:init:cond:p1}
and parameter values from Table~\ref{parameters}
with $\xi = 50$ and $T = 50$ (Section~\ref{sec:pb1:xineq0}).
The corresponding optimal control $u^*(t)$ is shown in Figure~\ref{fig:control}.}
\label{fig:spray}
\end{figure}

\begin{figure}[!htb]
\centering
\includegraphics[scale=0.49]{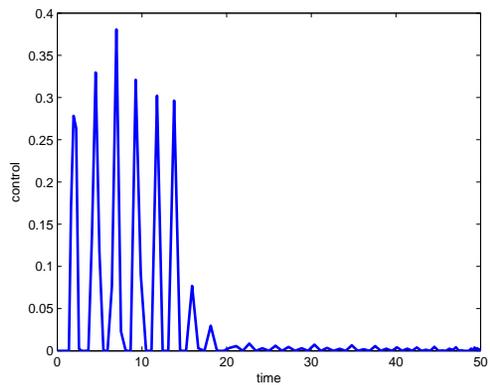}
\caption{The optimal spraying $u^*(t)$ for problem
\eqref{model:controls}--\eqref{eq:init:cond:p1}
and parameter values from Table~\ref{parameters}
with $\xi = 50$ and $T = 50$ (Section~\ref{sec:pb1:xineq0}).}
\label{fig:control}
\end{figure}

\begin{figure}[!htb]
\centering
\subfloat[$f(t)$]{\label{fig:preyswood2}
\includegraphics[scale=0.33]{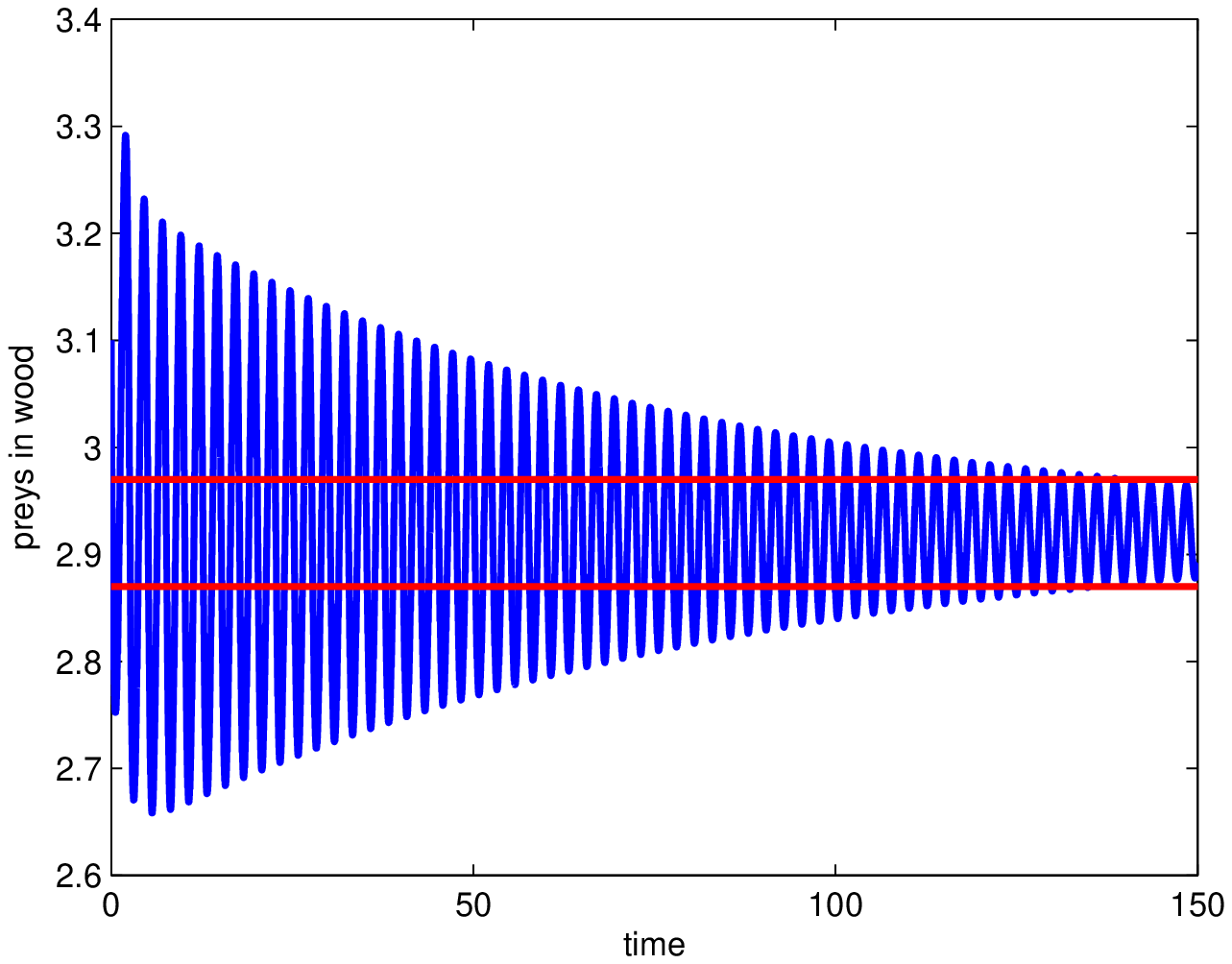}}
\subfloat[$s(t)$]{\label{fig:spiders2}
\includegraphics[scale=0.33]{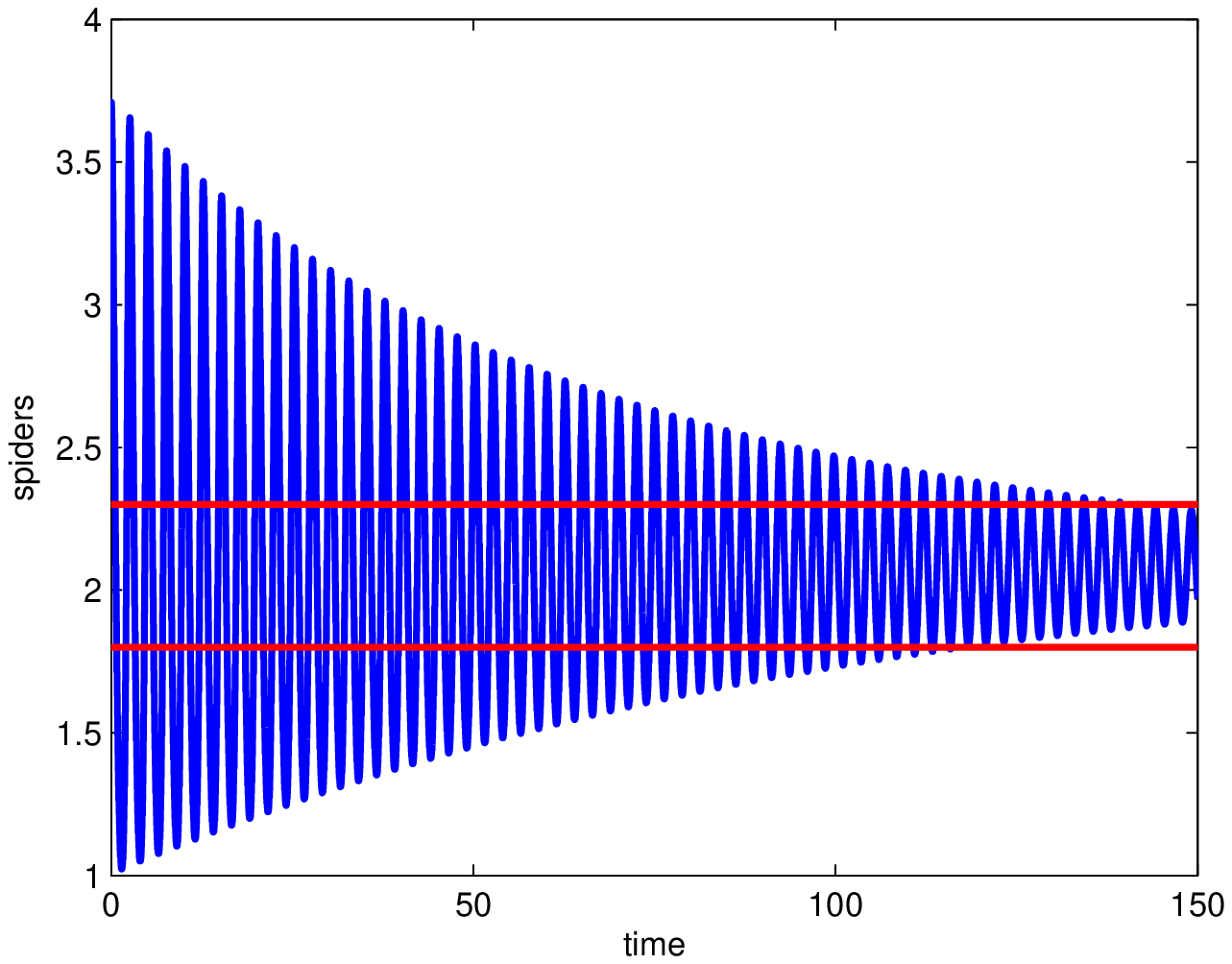}}
\subfloat[$v(t)$]{\label{fig:preyswine2}
\includegraphics[scale=0.33]{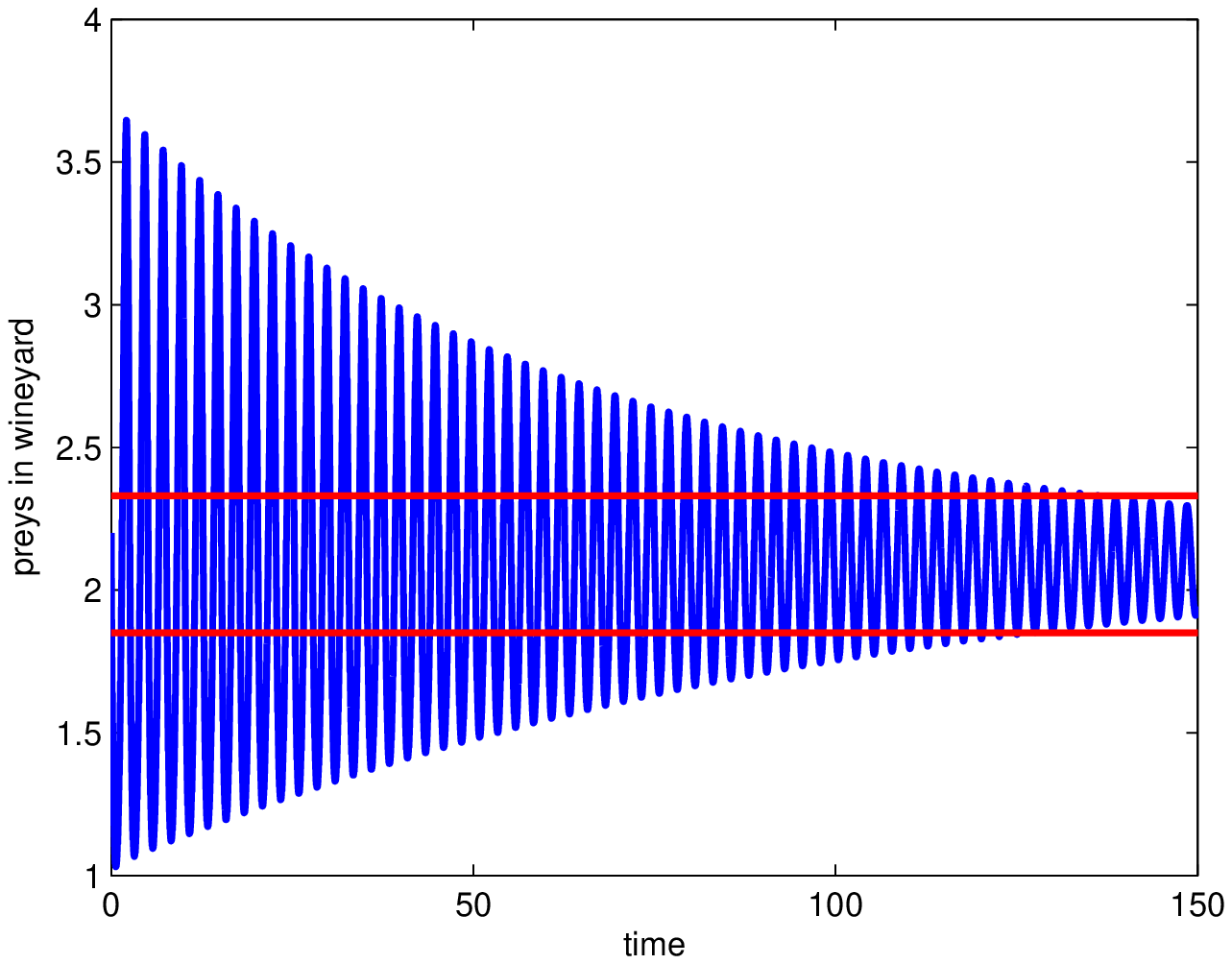}}
\caption{Behavior of \eqref{model:controls}
without spraying ($u(t) = 0$, $t \in [0, 150]$)
with initial conditions \eqref{eq:init:cond:p1} and 
parameter values from Table~\ref{parameters}.}
\label{fig:no:spray:range}
\end{figure}

Proceeding differently from Section~\ref{sec:xi0}, now we 
also account for the cost of insecticides ($\xi \ne 0$)
in the objective functional to be minimized.
More precisely, we set $\xi = 50$ in \eqref{costfunction}.
We observe that the coexistence equilibrium
is attained at $t \simeq 20$ (see Figure~\ref{fig:spray}).
The optimal way of spraying $u^*(t)$, $0 \le t \le 50$,
is illustrated in Figure~\ref{fig:control}.
Comparing the range of oscillations of $(f(\cdot), s(\cdot), v(\cdot))$
for $u \equiv 0$ (no spraying with insecticide) with those of
$(f^*(\cdot), s^*(\cdot), v^*(\cdot))$ corresponding to the optimal control $u^*$,
we observe that in case of human intervention
(the case with control $u^*$) the range of oscillations
at $t \simeq 26$ is similar to the case without control at $t \simeq 150$
(Figures~\ref{fig:spray} and \ref{fig:no:spray:range}).

In contrast with Section~\ref{sec:xi0}, now we take into account the costs
associated with the spray. For this reason the optimal control never attains
the maximum allowed value, being always less than $0.4$.
Moreover, for $t > 20$ days, the optimal control takes values close to zero,
which means that very small quantities of insecticides 
are applied (see Figure~\ref{fig:control}).
We observe that if no insecticide is applied,
one needs almost 150 days to attain the values
for $v(t)$ that we have at the end of 25 days
with insecticides (compare Figures~\ref{fig:spray}
and \ref{fig:no:spray:range}).


\subsection{A time-optimal control problem}
\label{sec:ocp2}

We now consider the objective functional in the form
\begin{equation}
\label{costfunction:mintf}
\mathcal{J} = \int_0^{T} 1 \, dt \, .
\end{equation}
Our aim is to minimize the time $T$ for which
the number of parasites living in a vineyard will vanish: $v(T) = 0$.
More precisely, the optimal control problem consists in finding the control $u^*$
that attains the final point ($f(T)=f_T$, $s(T)=s_T$, $v(T)=v_T$) in minimal time,
where $f_T$, $s_T > 0$ are free and $v_T=0$.
In this section, we make the following controllability assumption:

\medskip

$(H)$ The final point $\left(f_T, s_T, 0 \right)$, with $f_T$, $s_T > 0$,
can be steered from the initial point \eqref{init:cond}.

\begin{thrm}
\label{the:thm2}
Under assumption $(H)$, the optimal control problem
\eqref{model:controls}, \eqref{costfunction:mintf} admits a solution
$(f^*(\cdot), s^*(\cdot), v^*(\cdot))$ on $[0, T]$, associated to a control
$u^*(\cdot)$, where $T>0$ is the minimal time.
Moreover,
\begin{equation*}
u^*(t) = \min \left\{ \max \left\{0, u^*_{sing}(t) \right\}, 1 \right\},
\end{equation*}
where $u^*_{sing}$ denotes the singular control given by
\begin{equation*}
\begin{split}
u^*_{sing}(t) = -\frac{A(t)}{2 B(t)}
\end{split}
\end{equation*}
with
\begin{equation*}
\begin{split}
A&= q{b}^{2}{s}^{2}{W}^{2}{V}^{2}(\lambda_3-\lambda_2\,k)
+ {W}^{2}{V}^{2}q(e-bs) +2\,q\lambda_3\,ev{W}^{2}bsV+2\,
\lambda_1\,{r}^{2}{f}^{2}{V}^{2}+{W}^{2}{V}^{2}({r}^{2}\lambda_1\\
&\quad+2\,csq\lambda_1\,r+{c}^{2}{s}^{2}q\lambda_2\,k)+2\,csW{V}^{2}\lambda_1\,rf
+{W}^{2}{V}^{2}(-2\,q e\lambda_3\,bs+r\lambda_2\,skc+Kq{a}^{2}\lambda_2-\lambda_1\,csa\\
&\quad-{c}^{2}{s}^{2}q\lambda_1)+2qVW(\lambda_1\,{r}^{2}f V-\lambda_3\,{e}^{2}v W)
+{W}^{2}{V}^{2}(-2\,cs\lambda_1\,r-{c}^{2}{s}^{2}\lambda_2\,k-Kqbv)\\
&\quad-2(\lambda_1\,{r}^{2}{f}^{2}q{V}^{2}+ q\lambda_3\,{e}^{2}{v}^{2}{W}^{2})
+{W}^{2}{V}^{2}(-{r}^{2}q\lambda_1+{c}^{2}{s}^{2}\lambda_1+q{e}^{2}\lambda_3)\\
&\quad-2 qev{W}^{2}V-2\,\lambda_1 {r}^{2}f{V}^{2}W
-{W}^{2}{V}^{2}skc\lambda_3 bv+
{W}^{2}{V}^{2}Kqa(\lambda_1 cf-2\lambda_2 kbv
-2 \lambda_2 kcf + \lambda_3 bv)\\
&\quad+\lambda_1{W}^{2}{V}^{2}cskbv-2\,csqW{V}^{2}\lambda_1\,rf
+{W}^{2}{V}^{2}\lambda_2 s k(qeb-rqc)+2\lambda_2\,sk(rfq{V}^{2}cW-evq{W}^{2}bV \\
&\quad- rf{V}^{2}cW)+{W}^{2}{V}^{2}Kq(-kbv\lambda_1\,cf
+{k}^{2}{b}^{2}{v}^{2}\lambda_2+2
\,{k}^{2}bv\lambda_2\,cf- k{b}^{2}{v}^{2}\lambda_3
-k{c}^{2}{f}^{2}\lambda_1+{k}^{2}{c}^{2}{f}^{2}\lambda_2\\
&\quad-kcf\lambda_3\,bv +sc\lambda_3\,bv)
+{W}^{2}{V}^{2}\lambda_1(-qskb\,cf
+csqa-csqkbv)+q\lambda_3\,bs{W}^{2}{V}^{2}(-a+kcf)\\
&\quad+Kq(c{f}^{2}W{V}^{2}\lambda_1\,r+b{v}^{2}{W}^{2}V
\lambda_3\,e+\lambda_2\,k(-bv{W}^{2}{V}^{2}e+b{v}^{2}{W}^{2}Ve-cf{W}^{2}{V}^{2}r
+c{f}^{2}W{V}^{2}r)),
\end{split}
\end{equation*}
\begin{equation*}
\begin{split}
B&=hWV [\lambda_2 K WV kc q(1-q)+\lambda_1 V (2\,rq
-r{q}^{2}- cKqW + cK{q}^{2}W -r)\\
&\qquad +KWV{q}^{2}b(\lambda_2\,k -\lambda_3)-{q}^{2}\lambda_3\,eW ] \, ,
\end{split}
\end{equation*}
and
\begin{equation}
\label{eq:adj:vect:pmp:tfmin}
\begin{cases}
\dot{\lambda}_1(t) = -\lambda_1(t) \left( r \left( 1-{\frac {f^*(t)}{W}} \right)
-{\frac {rf^*(t)}{W}}-c s^*(t) \right) - \lambda_2(t) \,s^*(t) k c,\\
\dot{\lambda}_2(t) = \lambda_1(t)\,cf^*(t) - \lambda_2(t) \,
\left( -a + k b v^*(t)+k c f^*(t) \right) + \lambda_3(t) \,b v^*(t),\\
\dot{\lambda}_3(t) = - \lambda_2(t) \,s^*(t) k b- \lambda_3(t)
\, \left( e \left( 1-{\frac {v^*(t)}{V}} \right) -{\frac {e v^*(t)}{V}}-b s^*(t) \right).
\end{cases}
\end{equation}
\end{thrm}

\begin{proof}
The existence of a solution $(f^*(\cdot), s^*(\cdot), v^*(\cdot))$,
associated to a control $u^*(\cdot)$ on $[0, T]$, where $T>0$ is the minimal time,
comes from the assumption $(H)$ (see, e.g., \cite[Chapter 9]{Cesari_1983}
for optimal control existence theorems).
According to the Pontryagin Maximum Principle \cite{Pontryagin_et_all_1962},
there exist a real number $\lambda^0 \geq 0$ and an adjoint function
$\lambda(\cdot) = (\lambda_1(\cdot), \lambda_2(\cdot), \lambda_3(\cdot))
\, : \, [0, T] \rightarrow \mathbb{R}^3$, with
$\left( \lambda_1(\cdot), \lambda_2(\cdot),
\lambda_3(\cdot), \lambda^0 \right) \neq 0$,
such that \eqref{eq:adj:vect:pmp:tfmin} holds
together with the minimality condition
\begin{equation}
\label{mincondPMP:tfmin}
H(f^*(t), s^*(t), v^*(t), u^*(t), \lambda(t), \lambda^0)
= \min_{0 \leq u \leq 1} H(f^*(t), s^*(t), v^*(t), u, \lambda(t), \lambda^0)
\end{equation}
holding almost everywhere on $[0, T]$. Here
\begin{equation*}
\begin{split}
H &= H(f, s, v, u, \lambda, \lambda^0)\\
&= \lambda^0 +\lambda_1 \left( r f\left(1 - \frac{f}{W}\right)
- c s f - h(1-q) u \right)\\
&\quad +\lambda_2 \left( s\left(-a + k b v
+ k c f\right) - h K q u \right)\\
&\quad +\lambda_3 \left( e v\left(1-\frac{v}{V}\right)
- b s v - h q u \right) \, .
\end{split}
\end{equation*}
Moreover, $\min_{0 \leq u \leq 1} H(f^*(t), s^*(t), v^*(t), u, \lambda(t), \lambda^0) = 0$
for every $t \in [0, T]$. From \eqref{mincondPMP:tfmin} we have
\begin{equation*}
u^*(t) = \min \left\{ \max \left\{0, u^*_{sing} \right\}, 1 \right\},
\end{equation*}
where $u^*_{sing}$ denotes a singular control that can be obtained by differentiating twice
$\frac{\partial H}{\partial u} =0$, i.e., 
$$
h(1-q)-\lambda_2 h K q - \lambda_3 h q  = 0,
$$
and replacing each derivative $\dot{\lambda}_i$, $i =1,2,3$,
by \eqref{eq:adj:vect:pmp:tfmin}, and $(\dot{f}^*, \dot{s}^*, \dot{v}^*)$ by \eqref{model:controls}.
\end{proof}

\begin{figure}[!htb]
\centering
\includegraphics[scale=0.49]{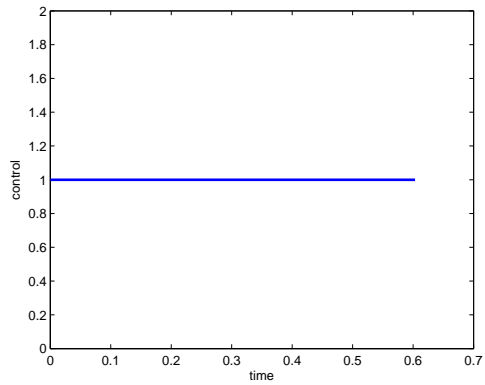}
\caption{The optimal spraying $u^*(t)$ for the minimal time problem of Section~\ref{sec:ocp2}.
Note that the spraying stops at time $T=0.6$ because at this time the vineyard
insect population vanishes, see Figure \ref{fig:spray:mintf}.}
\label{fig:control:mintf}
\end{figure}

\begin{figure}[!htb]
\centering
\subfloat[$f^*(t)$]{\label{fig:preyswood:op2}
\includegraphics[scale=0.33]{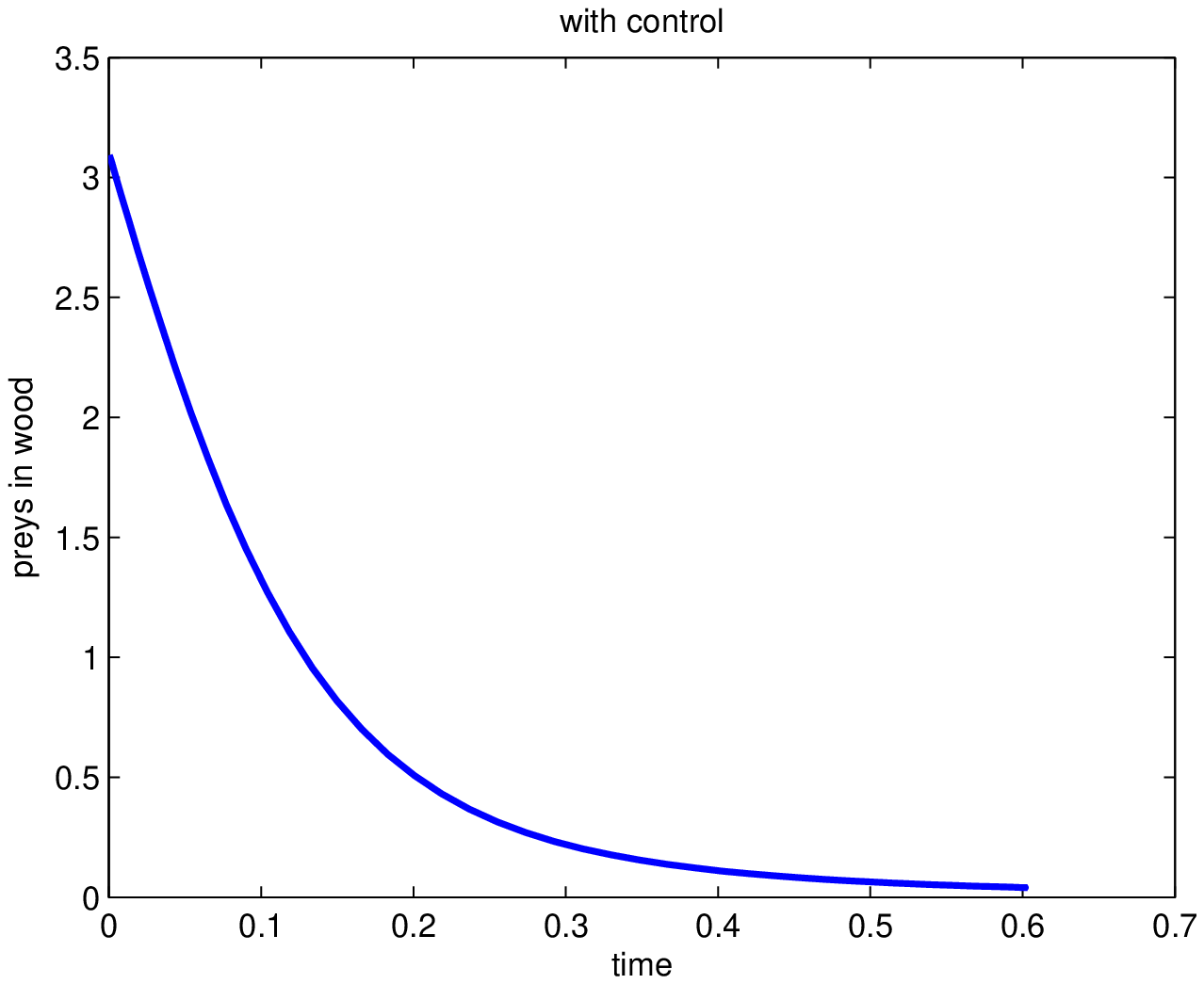}}
\subfloat[$s^*(t)$]{\label{fig:spiders:op2}
\includegraphics[scale=0.33]{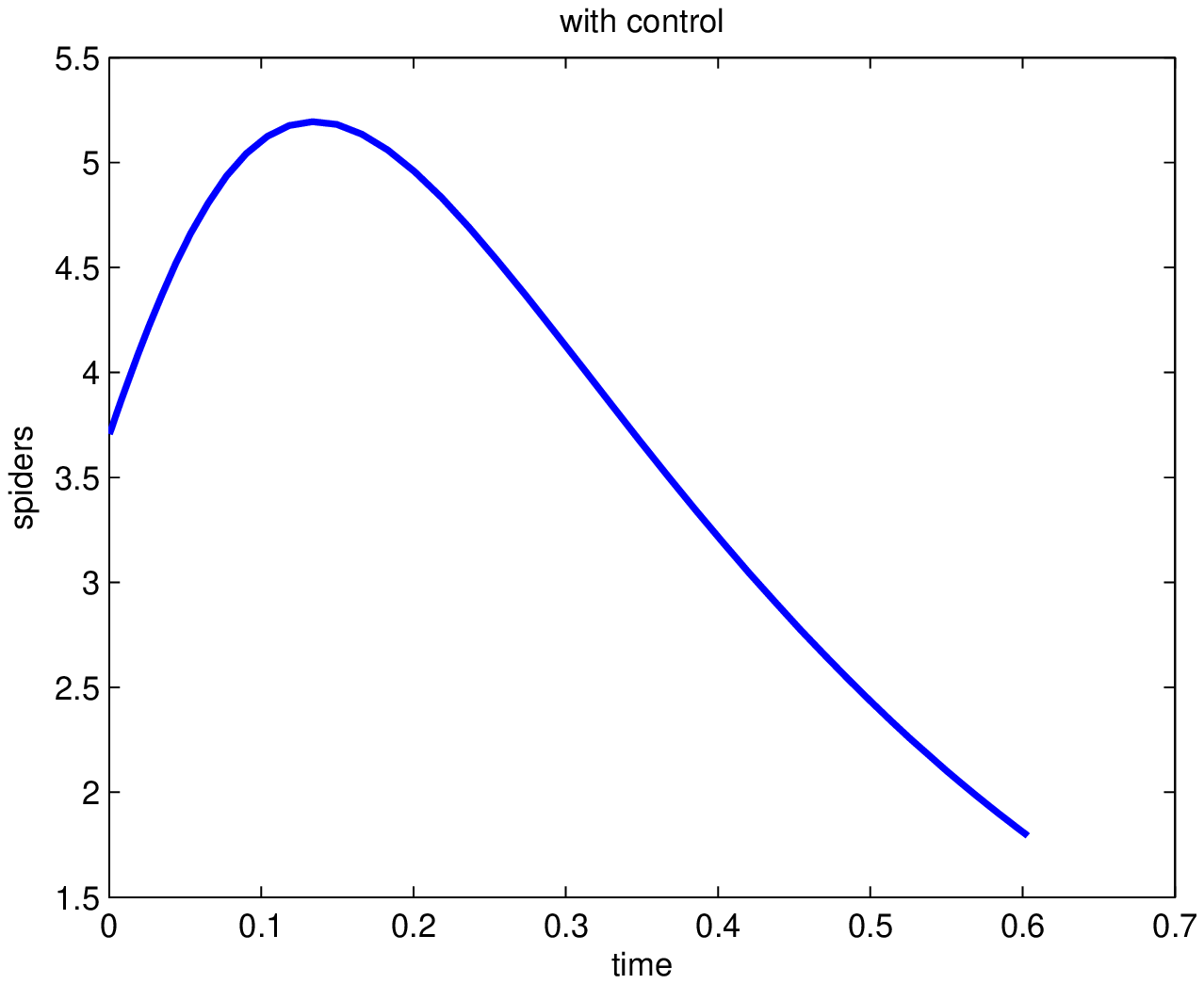}}
\subfloat[$v^*(t)$]{\label{fig:preyswine:op2}
\includegraphics[scale=0.33]{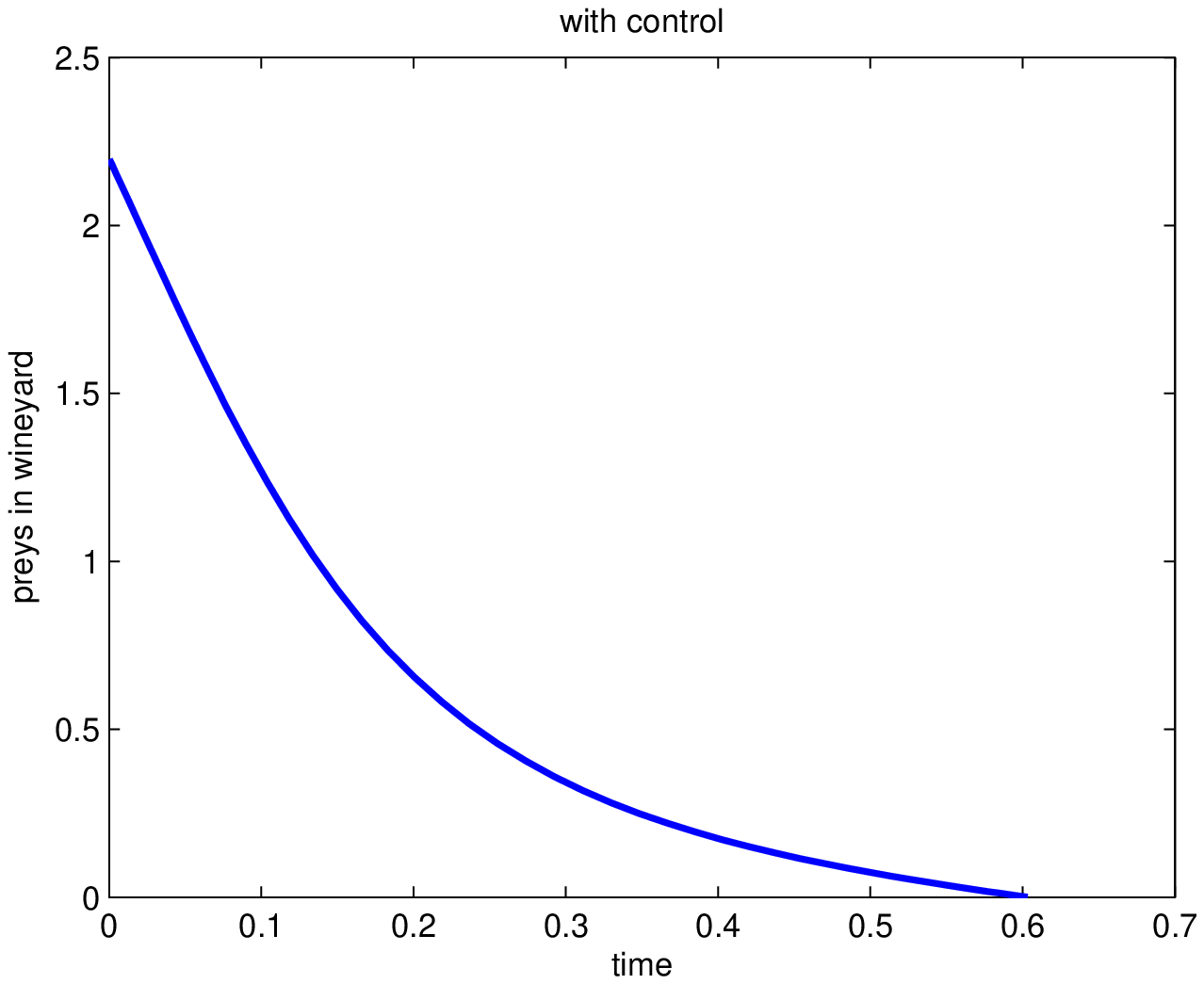}}
\caption{Behavior of system \eqref{model:controls} when subject to the
optimal control spraying $u^*(\cdot)$, solution
of the minimal time problem of Section~\ref{sec:ocp2}.}
\label{fig:spray:mintf}
\end{figure}

\begin{figure}[!htb]
\centering
\subfloat[$f(t)$]{\label{fig:preyswood:op3:a}
\includegraphics[scale=0.33]{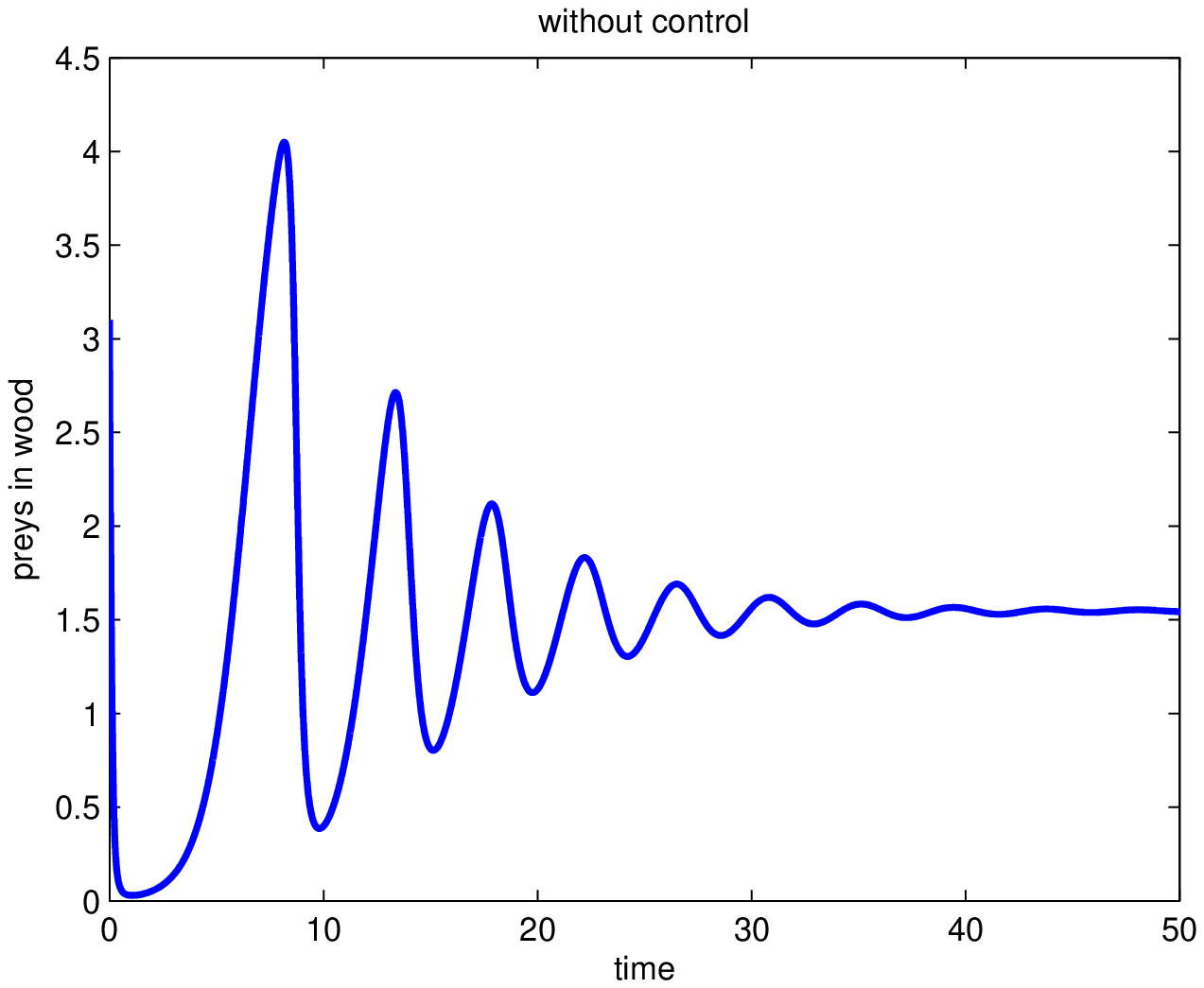}}
\subfloat[$s(t)$]{\label{fig:spiders:op3:a}
\includegraphics[scale=0.33]{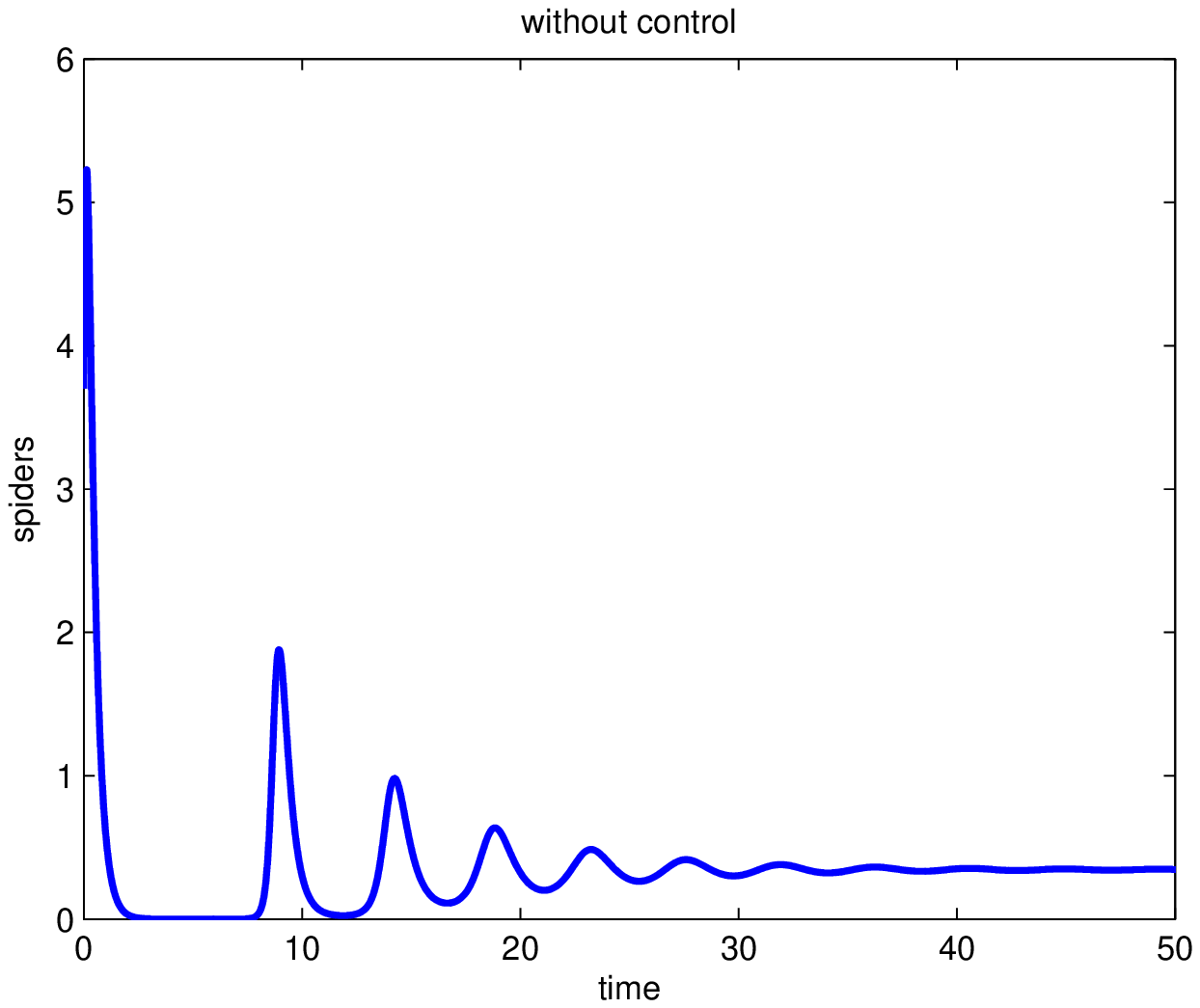}}
\subfloat[$v(t)$]{\label{fig:preyswine:op3:a}
\includegraphics[scale=0.33]{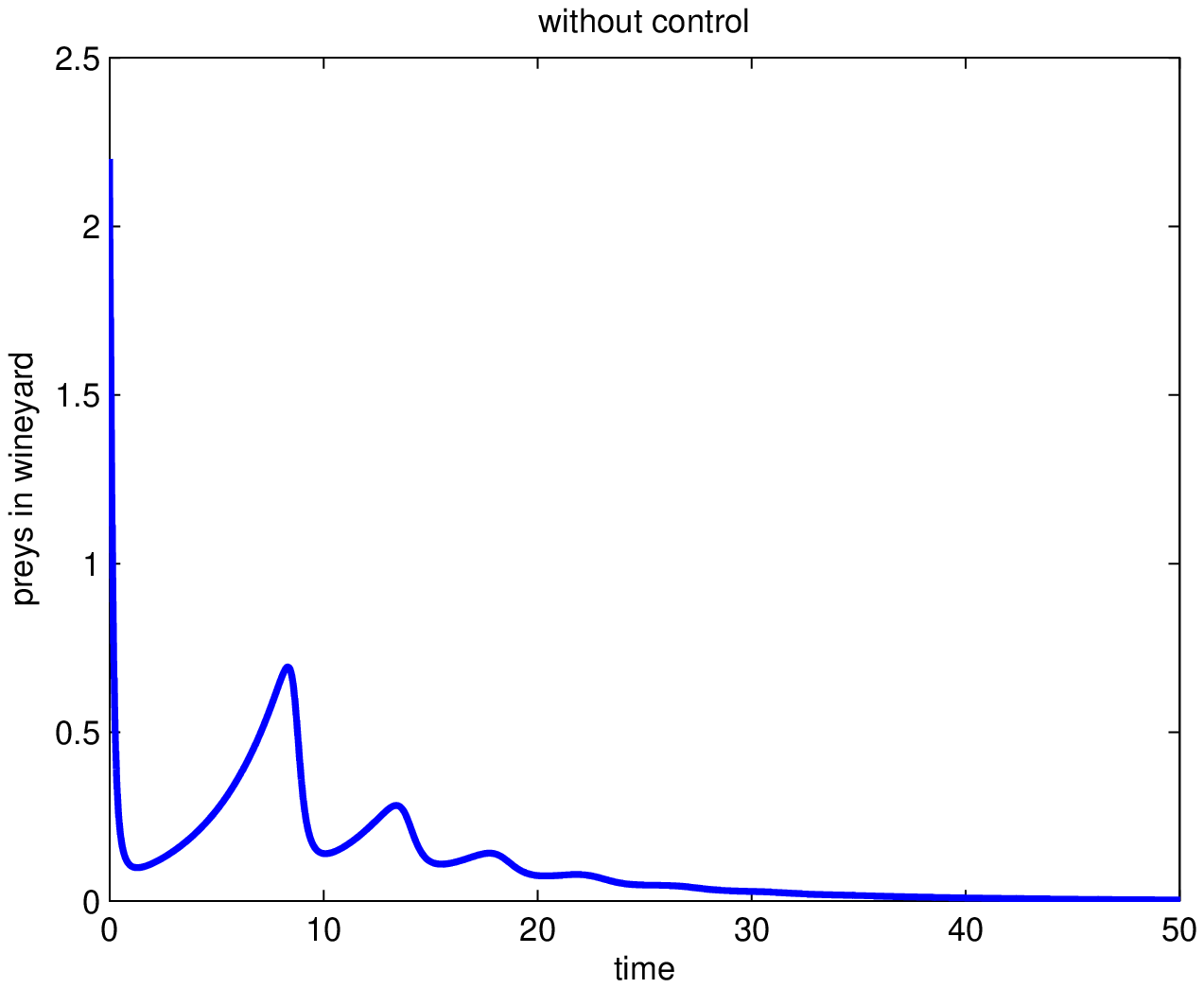}}
\caption{Behavior of system \eqref{model:controls}
without spraying ($u(t) = 0$, $t \in [0,50]$),
subject to \eqref{eq:init:cond:p1} and parameter values from
Table~\ref{parameters} with the exception of $e = 0.3$.}
\label{fig:sol:mintf}
\end{figure}

For numerical simulations, we considered $e = 0.3$ and the rest
of the parameter values from Table~\ref{parameters},
which satisfy the conditions for the feasibility and stability of the equilibrium
$(f_4, s_4, 0)$, that is,
$$
a < ckW \ \text{ and } \ c^2kWe < br(ckW - a)
$$
(see \cite[pag.~197]{Ezio:spiders:JNAIAM:2008}). As for the initial conditions,
we considered the same \eqref{eq:init:cond:p1} as in the previous section.
Figure~\ref{fig:control:mintf} shows the optimal control,
and Figure~\ref{fig:spray:mintf} the optimal state variables $(f^*, s^*, v^*)$.
In Figure~\ref{fig:sol:mintf} we can observe the solutions of model \eqref{model:controls}
without control, i.e., $u\equiv 0$ with initial conditions \eqref{eq:init:cond:p1}
and parameter values from Table~\ref{parameters} with the exception of $e = 0.3$.


\section{Discussion}
\label{sec:disc}

Comparing the no-control policy with the control for minimization of the total number
of pests, but not accounting for the costs of the pesticides,
we find that in the latter case both spiders and pests experience
outbursts during the last part of the simulation, higher than the peaks of the
oscillations arising naturally in the absence of the control, compare Figures
\ref{fig:no:spray} and \ref{fig:spray:xi0}.

On the same problem, we replicated the simulation over a larger timeframe, $T=150$.
In this case we note that the resulting control is different. Further, it allows
again population peaks toward the final time, but their sizes now are comparable to
the ones of the system without control, compare Figures \ref{fig:control:xi0}
and \ref{fig:no:spray}. For earlier times instead, the populations remain at lower values. 
This apparently indicates that the longer the timespan of the control application, 
the lower the possible populations outbursts are.

Taking into account the costs of spraying, it is immediately seen that all the populations
after a relatively initial short transient, up to $t=20$,
are now confined to small oscillations about the average value that they attain when
the control is not administered, compare Figures \ref{fig:spray}
and \ref{fig:no:spray}. Note also that the corresponding control is mainly 
used at the early times, and is essentially ``negligible'' after
time $t=20$, see Figure \ref{fig:control}.
Further, the optimal spraying policy uses all the time
only at most $40\%$ of the available insecticide. While in the previous cases,
Figure~\ref{fig:control:xi0} and, above all, Figure~\ref{fig:control:xi0:tf150}, 
the control is maximal and used constantly for a good portion of the time interval,
here it is administered only at a few selected instants, compare the peaks 
of Figure~\ref{fig:control} with the flat portion of Figure~\ref{fig:control:xi0:tf150}.
The same result is achieved in Figure~\ref{fig:no:spray:range}, where
no control is used, but it takes much more time, in fact populations fall within 
the same range only at time $t=150$, rather than $t=20$ as found above.

Finally, comparing Figures~\ref{fig:spray:mintf} and \ref{fig:sol:mintf},
we see that the minimal time under spraying is about a hundredfold smaller 
than without the insecticide usage. The optimal control,
Figure~\ref{fig:control:mintf}, needs to be used only for the first $0.6$ 
units of time, time at which the pests are eradicated,
after which the system can keep to evolve naturally.


\section{Conclusion and future work}
\label{sec:conc:fw}

In this paper, we have considered two optimal control problems for a model of
pest management in agroecosystems. More precisely, we studied  
problems of how to spray pests with the aim to eradicate the parasites.
We proved necessary optimality conditions 
and obtained numerical solutions for each of the problems
using the optimal control solver PROPT \cite{PROPT}.
The problems considered in Section~\ref{sec:2} 
are challenging, both from analytical and numerical points of view.
Indeed, in Section~\ref{sec:xi0}, we consider a problem formulation 
that falls into the category of $L^1$-optimization. 
Numerical algorithms are especially poor in locating the optimal solutions for 
$L^1$-type problems, which is the main reason for the abundance 
of papers written on optimization with $L^2$-type objectives \cite{MyID:353}. 
In this work we have used the state of the art on optimal control solvers,
provided by PROPT (see our code in Appendix~\ref{append:propt:code}).
It remains an open problem to prove that PROPT candidates  
for local minimum are indeed (local) minimizers. We are not aware
of any kind of solvers that guarantee minimality, either local or global,
for our specific problems. The difficulties are related with the occurrence
of singular arcs, which may not be located by numerical procedures. 

For future work, it would be interesting to validate the obtained
numerical results with the necessary optimality results 
given by Theorems~\ref{the:thm} and \ref{the:thm2}.
Moreover, to prove sufficient optimality conditions and/or
global minimality remain also nontrivial open questions. 


\begin{acknowledgement}
This work was partially supported by project TOCCATA, 
reference PTDC/EEI-AUT/2933/2014, funded by Project 
3599 -- Promover a Produ\c{c}\~ao Cient\'{\i}fica e Desenvolvimento
Tecnol\'ogico e a Constitui\c{c}\~ao de Redes Tem\'aticas (3599-PPCDT)
and FEDER funds through COMPETE 2020, Programa Operacional
Competitividade e Internacionaliza\c{c}\~ao (POCI), and by national
funds through Funda\c{c}\~ao para a Ci\^encia e a Tecnologia (FCT)
and CIDMA, within project UID/MAT/04106/2013 (Silva and Torres);
the post-doc fellowship SFRH/BPD/72061/2010 (Silva);
and by the project ``Metodi numerici in teoria delle popolazioni''
of the Dipartimento di Matematica ``Giuseppe Peano'' (Venturino).
The authors are grateful to two referees for useful comments and suggestions.
\end{acknowledgement}



\appendix


\section{PROPT Matlab codes}
\label{append:propt:code}

We present here our PROPT Matlab codes
so that any interested reader can replicate the numerical
results reported in this work. We begin with the code for
the optimal control problem investigated in Section~\ref{sub:sec:ocp1},
where the parameter values are given in Table~\ref{parameters}.

{\small
\begin{verbatim}
% Problem setup

toms t
tf = 150; 
p = tomPhase('p', t, 0, tf, 101);  
setPhase(p);

tomStates x1 x2 x3

tomControls u1

x = [x1; x2; x3];

u = [u1];

% Initial conditions 

x0i = [3.1; 3.7; 2.2];

x0 = icollocate({x1==x0i(1),x2==x0i(2),x3==x0i(3)});

% Box constraints and boundary

uL = [0]; uU = [1];

cbox = {0 <= icollocate(x1) <= 1000
0 <= icollocate(x2) <= 1000
0 <= icollocate(x3) <= 1000
0 <= icollocate(u1) <= 1 };

cbb = {collocate(uL <= u1 <= uU)

initial(x == x0i)};


% Parameters of the model

a = 3.1; b = 1.2; c= 0.2; e = 2.5; r = 1; W = 5; V = 1000; k = 1; 
h = 0.7; K = 0.01; q = 0.9; 

% ----- control system -------------------

% f = x1; s = x2; v = x3; 

ceq = collocate({
dot(x1) == r.*x1.*(1 - x1/W) - c*x2.*x1 - h*(1-q)*u1;
dot(x2) ==  x2.*(-a + k*b*x3 + k*c*x1) - h*K*q*u1;
dot(x3) == e*x3.*(1 - x3/V) - b*x2.*x3 - h*q*u1});

% ----- cost functional ------------------

objectiveJ1 = integrate(x3 + 50*u1.^2);

% objectiveJ1 = integrate(x3);


% ------ Solve the problem --------

options = struct;

options.name = 'system predator prey';

solution = ezsolve(objectiveJ1, {cbb, cbox, ceq}, x0, options);
\end{verbatim}}

For the time-optimal control problem of Section~\ref{sec:ocp2},
the following code was developed:
{\small 
\begin{verbatim}
toms t
toms t_f
p = tomPhase('p', t, 0, t_f, 50);
setPhase(p);

tomStates x1 x2 x3

tomControls u

% Box constraints

cbox = {0 <= t_f <= 50
0 <= collocate(u) <= 1};

% Boundary constraints

cbnd = {initial({x1 == 3.1; x2 == 3.7; x3 == 2.2}) 
final({x3 == 0})};  

% Parameters of the model

a = 3.1; b = 1.2; c= 2; e = 0.3; r = 1; W = 5; V = 1000; k = 1; 
h = 0.7; K = 0.01; q = 0.9; 

% ODEs and path constraints

ceq = collocate({dot(x1) == r.*x1.*(1 - x1/W) - c*x2.*x1 - h*(1-q)*u;
dot(x2) ==  x2.*(-a + k*b*x3 + k*c*x1) - h*K*q*u;
dot(x3) == e*x3.*(1 - x3/V) - b*x2.*x3 - h*q*u});

% Objective

objective = t_f;

% Solve the problem

options = struct;
options.name = 'minimum time system predator prey';
options.prilev = 1;

solution = ezsolve(objective, {cbox, cbnd, ceq}, options);
\end{verbatim}}



\begin{thebibliography}{99}

\bibitem{Behncke2000}
H. Behncke.
\textit{Optimal control of deterministic epidemics}.
Optimal Control, Applications and Methods, 21:269--285, 2000.

\bibitem{Cesari_1983}
L. Cesari.
\textit{Optimization---theory and applications}.
Applications of Mathematics (New York), 17,
Springer, New York, 1983.

\bibitem{Ezio:spiders:JNAIAM:2008}
S. Chatterjee, M. Isaia, F. Bona, G. Badino\ and\ E. Venturino.
\textit{Modelling environmental influences on wanderer
spiders in the Langhe region (Piemonte-NW Italy)}.
JNAIAM J. Numer. Anal. Ind. Appl. Math., 3:3-4, 193--209, 2008.

\bibitem{CIVa}
S. Chatterjee, M. Isaia\ and\ E. Venturino.
\textit{Effects of spiders' predational delays in intensive agroecosystems}.
Nonlinear Anal. Real World Appl., 10:5, 3045--3058, 2009.

\bibitem{CIVb}
S. Chatterjee, M. Isaia\ and\ E. Venturino.
\textit{Spiders as biological controllers in the agroecosystem}.
J. Theoret. Biol., 258:3, 352--362, 2009.

\bibitem{CD98}
M. J. Costello\ and\ K. M. Daane.
\textit{Influence of ground cover on spider populations in a table grape vineyard}.
Ecological Entomology, 23:1, 33--40, 1998.

\bibitem{RD-B}
P. DeBach\ and\ D. Rosen.
\textit{Biological control by natural enemies}.
Cambridge Univ. Press, 2nd ed., Cambridge, 1991.

\bibitem{Eisen1979}
M. Eisen.
\textit{Mathematical Models in Cell Biology and Cancer Chemotherapy}.
Lectures Notes in Biomathematics, Vol. 30, Springer Verlag, 1979.

\bibitem{Fleming_Rishel_1975}
W. H. Fleming\ and\ R. W. Rishel.
\textit{Deterministic and stochastic optimal control}.
Springer, Berlin, 1975.

\bibitem{Gaff2009}
H. Gaff\ and\ E. Schaefer.
\textit{Optimal control applied to vaccination and 
treatment strategies for various epidemiologic models}.
Math. Biosci. Eng., 6:3, 469--492, 2009.

\bibitem{Winemaking}
E. M. Harkness, R. P. Vine\ and\ S. J. Linton.
\textit{Winemaking: From Grape Growing to Marketplace}.
Springer, 2002.

\bibitem{HattafHIV2012} 
K. Hattaf\ and\ N. Yousfi.
\textit{Two optimal treatments of HIV infection model}. 
World Journal of Modelling and Simulation, 8:27--35, 2012.

\bibitem{HattafHIVdelay2012} 
K. Hattaf\ and\ N. Yousfi.
\textit{Optimal control of a delayed HIV infection model 
with immune response using an efficient numerical method}. 
ISRN Biomathematics, Art. ID 215124, 7~pp, 2012.

\bibitem{IbrahimIJDC2016}
F. Ibrahim, K. Hattaf, F. A. Rihan\ and\ S. Turek.
\textit{Numerical method based on extended one-step schemes 
for optimal control problem with time-lags}.
Int. J. Dynam. Control, in press,
doi:10.1007/s40435-016-0270-x 

\bibitem{IBBB}
M. Isaia, G. Badino, F. Bona\ and\ E. Bosca.
\textit{I ragni costruttori di tela nella valutazione della qualit\`a ambientale:
un esempio di applicazione, (Webbuilder spiders in evaluating the environmental quality:
example of an application)}. 
In: Atti XXIII convegno S.It.E., Ecologia quantitativa,
8-10 Sept. 2003, Como, Italy, 
Edizione speciale Vincitori Premio Marchetti, 27:61--66, 2004.

\bibitem{IBB06}
M. Isaia, F. Bona\ and\ G. Badino.
\textit{Influence of landscape diversity and agricultural practices 
on spiders assemblage in Italian vineyards of Langa Astigiana (Northwest Italy)}.
Environ. Entomol., 35:2, 297--307, 2006.

\bibitem{SLenhart_2002}
E.~Jung, S.~Lenhart\ and\ Z.~Feng.
\textit{Optimal control of treatments in a two-strain tuberculosis model}.
Discrete Contin. Dyn. Syst. Ser. B, 2:4, 473--482, 2002.

\bibitem{Ledzewicz2002}
U. Ledzewicz\ and\ H. Sch\"attler.
\textit{Optimal bang-bang controls for a two-compartment model in cancer chemotherapy}.
J. Optim. Theory Appl., 114:3, 609--637, 2002.

\bibitem{Lenhart2007}
S. Lenhart\ and\ J. T. Workman.
\textit{Optimal control applied to biological models}.
Chapman \& Hall/CRC, Boca Raton, FL, 2007.

\bibitem{Ma}
P. Marc.
\textit{Intraspecific predation in Clubiona corticalis (Walckenaer, 1802) 
(Araneae, Clubionidae): a spider bred for its interest in biological control}.
Mem. Queensl. Mus., 33:2, 607--614, 1993.

\bibitem{MC}
P. Marc\ and\ A. Canard.
\textit{Maintaining spider biodiversity in agroecosystems as a tool in pest control}.
Agric. Ecosyst. Environ., 62:2-3, 229--235, 1997.

\bibitem{NM}
T. Nobre\ and\ C. Meierrose.
\textit{The species composition, within-plant distribution, and possible predatory
role of spiders (Araneae) in a vineyard in southern Portugal}.
In: Proceedings of the 18th European Colloquium on Arachnology
(Eds. P. Gajdos and S. Pek\^ er), Ekol\'{o}gia (Bratislava), 19:3, 193--200, 2000.

\bibitem{Pinney:HW}
T. Pinney.
\textit{A History of Wine in America: From the Beginnings to Prohibition}.
University of California Press, Berkeley, 1989.

\bibitem{Pontryagin_et_all_1962}
L. S. Pontryagin, V. G. Boltyanskii, R. V. Gamkrelidze\ and\ E. F. Mishchenko.
\textit{The mathematical theory of optimal processes}.
Translated from the Russian by K. N. Trirogoff; edited by L. W. Neustadt,
Interscience Publishers John Wiley \& Sons, Inc.\, New York, 1962.

\bibitem{RB}
S. E. Riechert\ and\ L. Bishop.
\textit{Prey control by an assemblage of generalist predators: 
spiders in garden test systems}.
Ecology, 71:4, 1441--1450, 1990.

\bibitem{RL}
S. E. Riechert\ and\ T. Lockley.
\textit{Spiders as biological control agents}.
Ann. Rev. Entomol., 29:299--320, 1984.

\bibitem{Ox:Comp:Wine}
J. Robinson.
\textit{The Oxford Companion to Wine}.
Oxford University Press, Oxford, 2006.

\bibitem{Sofia2010}
H. S. Rodrigues, M. T. T. Monteiro\ and\ D. F. M. Torres.
\textit{Dynamics of dengue epidemics when using optimal control}.
Math. Comput. Modelling, 52:9-10, 1667--1673, 2010.
{\tt arXiv:1006.4392}

\bibitem{Sofia2013}
H. S. Rodrigues, M. T. T. Monteiro\ and\ D. F. M. Torres.
\textit{Bioeconomic perspectives to an optimal control dengue model}.
Int. J. Comput. Math., 90:10, 2126--2136, 2013.
{\tt arXiv:1303.6904}

\bibitem{PaulaSilvaTorresT2014}
P. Rodrigues, C. J. Silva and D. F. M. Torres.
\textit{Cost-effectiveness analysis of optimal control measures for tuberculosis}.
Bull. Math. Biol., 76:10, 2627--2645, 2014.
{\tt arXiv:1409.3496}

\bibitem{PROPT}
P. E. Rutquist\ and\ M. M. Edvall.
\textit{PROPT -- Matlab Optimal Control Software}.
Tomlab Optimization, 2010. 

\bibitem{MyID:353}
C. J. Silva, H. Maurer\ and\ D. F. M. Torres.
\textit{Optimal control of a tuberculosis model with state and control delays}.
Math. Biosci. Eng., 14:1, 321--337, 2017.
{\tt arXiv:1606.08721}

\bibitem{SilvaTorresNACO}
C. J. Silva\ and\ D. F. M. Torres.
\textit{Optimal control strategies for tuberculosis treatment: a case study in Angola}.
Numer. Algebra Control Optim., 2:3, 601--617, 2012.
{\tt arXiv:1203.3255}

\bibitem{SilvaTorresMBS2013}
C. J. Silva\ and\ D. F. M. Torres.
\textit{Optimal control for a tuberculosis model 
with reinfection and post-exposure interventions}.
Math. Biosci., 244:2, 154--164, 2013.
{\tt arXiv:1305.2145}

\bibitem{Silva:Torres:DCDS-A:2015}
C. J. Silva\ and\ D. F. M. Torres.
\textit{A TB-HIV/AIDS coinfection model and optimal control treatment}.
Discrete Contin. Dyn. Syst., 35:9, 4639--4663, 2015.
{\tt arXiv:1501.03322}

\bibitem{VIBCB}
E. Venturino, M. Isaia, F. Bona, S. Chatterjee\ and\ G. Badino.
\textit{Biological controls of intensive agroecosystems: 
Wanderer spiders in the Langa Astigiana}.
Ecological Complexity, 5:2, 157--164, 2008.

\bibitem{VIBITB}
E. Venturino, M. Isaia, F. Bona, E. Issoglio, V. Triolo\ and\ G. Badino.
\textit{Modelling the spiders ballooning effect on the vineyard ecology}.
Math. Model. Nat. Phenom., 1:1, 137--159, 2006.

\bibitem{W}
D. H. Wise.
\textit{Spiders in Ecological Webs}.
Cambridge University Press, Cambridge, 1993.

\end{thebibliography}
\end{document}